\documentclass[a4paper, 12pt, oneside]{article} 
\usepackage{graphicx} 
\usepackage{amsmath} 
\usepackage{amssymb} 
\usepackage{amsfonts} 
\usepackage{amsthm} 
\usepackage{longtable} 
\usepackage{multirow}
\usepackage{color}
\usepackage{makeidx}

\makeindex
\title{Regular Polygonal Complexes in Space, I
\footnote{Version of June 3, 2009}}
\author{Daniel Pellicer\\
York University\\
Toronto, Ontario, Canada M3J 1P3\\[.05in]
{\small and} \\[.05in]
Egon Schulte\thanks{Supported by NSA-grant 
H98230-07-1-0005}\\
Northeastern University\\
Boston, MA 02115, USA}
\date{ }

\newtheorem{lemma}{Lemma}[section] 
\newtheorem{theorem}[lemma]{Theorem}

\newtheorem{definition}[lemma]{Definition}

\renewcommand{\Gamma}{\varGamma} 
\renewcommand{\epsilon}{\varepsilon} 
 
\renewcommand{\hat}{\widehat}

\renewcommand{\geq}{\geqslant}

\newcommand{\K}{\mathcal{K}} 
\newcommand{\E}{\mathbb{E}^3}
\def\apeir{\mathop{\rm apeir}}

\begin{document}

\maketitle

\begin{abstract}
\noindent
A polygonal complex in euclidean $3$-space $\E$ is a discrete polyhedron-like structure with finite or infinite polygons as faces and finite graphs as vertex-figures, such that a fixed number $r\geq 2$ of faces surround each edge. It is said to be regular if its symmetry group is transitive on the flags. The present paper and its successor describe a complete classification of regular polygonal complexes in $\E$. In particular, the present paper establishes basic structure results for the symmetry groups, discusses geometric and algebraic aspects of operations on their generators, characterizes the complexes with face mirrors as the $2$-skeletons of the regular $4$-apeirotopes in $\E$, and fully enumerates the simply flag-transitive complexes with mirror vector $(1,2)$. The second paper will complete the enumeration.
\end{abstract}


\section{Introduction}
The study of highly symmetric polyhedra-like structures in ordinary euclidean $3$-space $\E$ has a long and fascinating history tracing back to the early days of geometry. With the passage of time, various notions of polyhedra have attracted attention and have brought to light new exciting classes of regular polyhedra including well-known objects such as the Platonic solids, Kepler-Poinsot polyhedra, Petrie-Coxeter polyhedra, or the more recently discovered Gr\"unbaum-Dress~polyhedra (see \cite{coxeter,crsp,d1,d2,gr1}).  

The radically new {\em skeletal\/} approach to polyhedra pioneered in \cite{gr1} is essentially graph-theoretical and has had an enormous impact on the field. Since then, there has been a lot of activity in this area, beginning with the full enumeration of the $48$ ``new" regular polyhedra in $\E$ by Gr\"unbaum~\cite{gr1} and Dress~\cite{d1,d2} (with a simpler approach described in the joint works \cite{ordinary,rap} of the second author with McMullen), moving to the full enumeration of chiral polyhedra in $\E$ in \cite{chiral1,chiral2} (see also \cite{pelwei}), then continuing with corresponding enumerations of regular polyhedra, polytopes, or apeirotopes (infinite polytopes) in higher-dimensional euclidean spaces by McMullen~\cite{pm,pm1,pm2} (see also \cite{ar,bra}). For a survey, see \cite{ms3}.

The present paper and its successor~\cite{pelsch} describe a complete classification of regular  polygonal complexes in euclidean space $\E$.  Polygonal complexes are discrete polyhedra-like structures made up from convex or non-convex, planar or skew, finite or infinite (helical or zig-zag) polygonal faces, always with finite graphs as vertex-figures, such that every edge lies in at least two, but generally $r \ge 2$ faces, with $r$ not depending on the edge. As combinatorial objects they are incidence complexes of rank $3$ with polygons as $2$-faces (see \cite{kom1,kom2}). Polyhedra are precisely the polygonal complexes with $r=2$. A polygonal complex is said to be {\em regular\/} if its full euclidean symmetry group is transitive on the flags. The two papers are part of a continuing program to classify discrete polyhedra-like structures by transitivity properties of their symmetry group. Characteristic of this program is the interplay of the abstract, purely combinatorial, aspect and the geometric one of realizations of complexes and their symmetries.

The paper is organized as follows. In Sections~\ref{definitions} and \ref{thegroup}, respectively, we begin with the notion of a regular polygonal complex $\K$ and establish basic structure results for its symmetry group $G(\K)$ in terms of generators for distinguished subgroups. It is found that either $\K$ is simply flag-transitive, meaning that its {\em full\/} symmetry group $G(\K)$ has trivial flag stabilizers, or that $\K$ has face-mirrors, meaning that $\K$ has planar faces lying in reflection mirrors and that $G(\K)$ has flag stabilizers of order $2$ (containing the reflection in the plane spanned by the face in a flag). Then in Section~\ref{facemir} we characterize the regular complexes with face mirrors as the $2$-skeletons of the regular $4$-apeirotopes in $\E$. There are exactly four such complexes, up to similarity. In Section~\ref{opera} we discuss operations on the generators of $G(\K)$ that allow us to construct new complexes from old and that help reduce the number of cases to be considered in the classification. Finally, in Sections~\ref{finmir12} and \ref{infinmir12}, respectively, we enumerate the simply flag-transitive regular polygonal complexes with finite or infinite faces and mirror vector $(1,2)$. Apart from polyhedra, there are exactly eight polygonal complexes of this kind, up to similarity. 

In the subsequent paper~\cite{pelsch} we investigate other possible mirror vectors and employ the operations to complete the classification.

\section{Polygonal complexes}
\label{definitions}

Geometric realizations of abstract polytopes or complexes have attracted a lot of attention. One specific approach is to start with an abstract object and then study ways of realizing its combinatorial structure in a more explicit geometric setting (for example, see \cite[Ch.5]{rap} and \cite{monw5}). In this paper, however, we take a more direct geometric approach and explicitly define the objects as polyhedra-like structures in an ambient space. In fact, throughout, this space will be euclidean $3$-space $\E$. 

Informally, a polygonal complex will consist of a family of vertices, edges and (finite or infinite) polygons, all fitting together in a way characteristic for geometric polyhedra or polyhedral complexes.

We say that a {\em finite polygon\/}, or simply {\em $n$-gon\/}, $(v_1, v_2, \dots, v_n)$ in euclidean $3$-space $\E$ is a figure formed by distinct points $v_1, \dots, v_n$, together with the line segments $(v_i, v_{i+1})$, for $i = 1, \dots, n-1$, and $(v_n, v_1)$. Similarly, an {\em infinite polygon\/} consists of a sequence of distinct points $(\dots, v_{-2},v_{-1}, v_0, v_1, v_2\dots)$ and of line segments $(v_i, v_{i+1})$ for each $i$, such that each compact subset of $\E$ meets only finitely many line segments. In either case we refer to the points as {\em vertices\/} and to the line segments as {\em edges\/} of the polygon.

\begin{definition}
\label{compdef}
A {\em polygonal complex}, or simply {\em complex}, $\K$ in $\E$ consists of a set $\cal V$ of points, called {\em vertices}, a set $\cal E$ of line segments, called {\em edges}, and a set $\cal F$ of polygons, called {\em faces}, such that the following properties are satisfied.\\[-0.3in]
\begin{itemize}\label{defcomplex}
\item[a)] The graph defined by $\cal V$ and $\cal E$, called the {\em edge graph\/} of $\K$, is connected.\\[-0.3in]
\item[b)] The vertex-figure of $\K$ at each vertex of $\K$ is connected. By the {\em vertex-figure\/} of $\K$ at a vertex $v$ we mean the graph, possibly with multiple edges, whose vertices are the neighbors of $v$ in the edge graph of $\K$ and whose edges are the line segments $(u,w)$, where $(u, v)$ and $(v, w)$ are edges of a common face of $\K$.\\[-0.3in]
\item[c)] Each edge of $\K$ is contained in exactly $r$ faces of $\K$, for a fixed number $r \geq 2$.\\[-0.3in]
\item[d)] $\K$ is {\em discrete\/}, meaning that each compact subset of $\E$ meets only finitely many faces of $\K$.\\[-.3in]
\end{itemize}
\end{definition}

We call a polygonal complex a (geometric) {\em polyhedron\/} if $r=2$. For general properties of finite or infinite polyhedra in $\E$ we refer to \cite[Ch.7E]{rap} and \cite{gr1,ordinary}. Note that the underlying combinatorial ``complex" given by the vertices, edges and faces of a polygonal complex $\K$, ordered by inclusion, is an incidence complex of rank $3$ in the sense of \cite{kom1,kom2}. (When dealing with incidence complexes of rank $3$, we are suppressing their improper faces:\  the unique minimum face of rank $-1$ and the unique maximum face of rank $3$.) In the next section we shall require some basic structure results about the automorphism groups of incidence complexes obtained in \cite{kom2}.

Note that the discreteness assumption in Definition~\ref{compdef}(d) implies that the vertex-figure at every vertex is finite, and that, for complexes with regular polygons as faces (with ``regular" to be defined below), each compact subset of $\E$ meets only finitely many vertices and edges (in fact, under the regularity assumption on the faces, this could have served as an alternative definition for discreteness). 

The (geometric) {\em symmetry group} $G = G(\K)$ of a polygonal complex $\K$ consists of all isometries of the affine hull of $\K$ that map $\K$ to itself. Note that this affine hull is $\E$, except when $\K$ is planar. The symmetry group $G(\K)$ can be viewed as a subgroup of the (combinatorial) {\em automorphism group} $\Gamma(\K)$ of $\K$, which consists of the combinatorial automorphisms of the underlying combinatorial (incidence) complex, that is, of the incidence preserving bijections of the set of vertices, edges and faces of $\K$. Throughout we use the term {\em group of $\K$} to mean the (full) symmetry group $G(\K)$ of $\K$. 

A polygonal complex $\K$ is called (geometrically) {\em regular} if its symmetry group $G(\K)$ is transitive on the flags of $\K$, or {\em combinatorially regular} if its automorphism group $\Gamma(\K)$ is transitive on the flags of $\K$. Recall here that a {\em flag} of $\K$ is a $3$-element set consisting of a vertex, an edge and a face of $\K$, all mutually incident. Two flags of $\K$ are said to be {\em $i$-adjacent\/} if they differ precisely in their element of rank $i$, that is, their vertex, edge, or face if $i=0$, $1$ or $2$, respectively. Every flag of a polygonal complex $\K$ is $i$-adjacent to exactly one other flag if $i=0$ or $1$, or exactly $r-1$ other flags if $i=2$.

The (symmetry) group $G(\K)$ of a regular polygonal complex $\K$ is transitive, separately, on the vertices, edges, and faces of $\K$. In particular, the faces are necessarily regular polygons, either finite, planar (convex or star-) polygons or non-planar ({\em skew\/}) polygons, or infinite, planar zigzags or helical polygons (see \cite[Ch.1]{crcp} or \cite{gr1,ordinary}). Note that a regular complex $\K$ cannot have faces which are linear apeirogons; in fact, by the connectedness assumptions on $\K$, this would force $\K$ to ``collapse" onto a single apeirogon.

Notice also that we explicitly allow the vertex-figures of $\K$ to have multiple edges, to account for the possibility that two adjacent edges of a face are adjacent edges of more than one face. On the other hand, if $\K$ is regular, then all edges of vertex-figures have the same multiplicity and all are single or double edges. In fact, if two adjacent edges of a face are adjacent edges of another face, then the symmetry of $\K$ that fixes these two edges and maps the first face to the second face must necessarily be the reflection in the plane spanned by these two edges. This also shows that regular complexes with planar faces have vertex-figures with single edges; that is, they are simple graphs. 

In practice we can think of the vertex-figure of $\K$ at a vertex $v$ as the graph whose vertices are given by the neighboring vertices of $v$ in the edge graph of $\K$ and whose edges are represented by straight line segments, with double edges (if any) lying on top of each other.

A priori there seems to be no tractable duality theory for polygonal complexes in space. At the combinatorial level there are no problems:\  each incidence complex has a dual which is also incidence complex, of the same rank as the original (see \cite{kom1,kom2}). The dual of an incidence complex associated with a polygonal complex in space has $r$ ``vertices" on an ``edge", such that the ``vertex-figures" are isomorphic to polygons. Configurations of points and lines in space (with $r$ points on every line) come to mind as candidates for duals, but there seems to be no obvious way of relating configurations and polygonal complexes geometrically. Moreover, as a reminder of the limitations of the duality concept, even the class of geometric polyhedra (that is when $r=2$) is not closed under geometric duality (see \cite[7E]{rap}):\  there are examples of regular polyhedra which do not have a geometric dual which is also a regular polyhedron (any helix-faced polyhedron gives an example).

The underlying edge graphs of regular complexes are examples of highly symmetric ``nets" as studied in crystal chemistry (for example, see \cite{del,okee,okhy}, and note in particular that \cite{okee} describes the nets of many regular polyhedra using the notation and nomenclature of \cite[Ch. 7E]{rap}). Nets are infinite, periodic geometric graphs in space that represent crystal structures, with vertices corresponding to atoms and edges to bonds. From a chemist's perspective, the arrangement of the atoms in space is a central piece of information about a chemical compound. The famous {\em diamond net\/} is the underlying net of the diamond form of carbon and of several other compounds. In the diamond net, each vertex $v$ is joined to exactly four neighboring vertices such that any two neighbors, along with $v$, are vertices of exactly two hexagonal ``rings" (fundamental circuits) of the net; in \cite[p. 299]{okhy}, the local configuration is described by the ``long Schl\"afli symbol" $6_{2}\cdot 6_{2}\cdot 6_{2}\cdot 6_{2}\cdot 6_{2}\cdot 6_{2}$ (there are six choices for the two neighbors and each contributes $2$ rings with $6$ vertices). The diamond net with its hexagonal ring structure will occur in our enumeration as the underlying edge-graph of the complex $\K_{7}(1,2)$ of (\ref{7gen12}). Its faces are skew hexagons, six around each edge, such that the vertex-figure is the double edge-graph of the tetrahedron (the six edges of the tetrahedral vertex-figure, each counted twice, contribute the six terms $6_2$ in the above symbol).

\section{The symmetry group}
\label{thegroup}

In this section we establish some basic structure theorems about the (symmetry) group $G(\K)$ of a regular (polygonal) complex $\K$.

Let $\K$ be a regular complex, and let $G:=G(\K)$ be its group. Let $\Phi := \{F_0, F_1, F_2\}$ be a fixed, {\em base}, flag of $\K$, where $F_0$ is a vertex, $F_1$ an edge and $F_2$ a face of $\K$. We denote the stabilizers in $G$ of $\Phi \setminus \{F_i\}$ and more generally of a subset $\Psi$ of $\Phi$ by
\[G_i = G_i(\K) := \{R \in G \,|\, F_j R = F_j \; \mbox{for $j \ne i$}\} \quad\; (i=0,1,2)\]
and
\[G_{\Psi} := \{R \in G \,|\, FR = F \; \mbox{for $F \in \Psi$}\},\]
respectively. Then $G_i = G_{\{F_j, F_k\}}$ with $i, j, k$ distinct, and $G_\Phi$ is the stabilizer of the base flag $\Phi$. Moreover,
\[G_\Phi = G_0 \cap G_1 = G_0 \cap G_2 = G_1 \cap G_2.\]
We also set $G_{F_i}:=G_{\{F_i\}}$ for $i=0,1,2$; this is the stabilizer of $F_i$. The stabilizer subgroups satisfy the following properties.

\begin{lemma}\label{groupstruct}
Let $\K$ be a regular complex with group $G=G(\K)$. Then, \\[-0.3in]
\begin{itemize}
\item[a)] $G = \langle G_0, G_1, G_2 \rangle$;\\[-.3in]
\item[b)] $G_\Psi = \langle G_j \,|\, F_j \notin \Psi\rangle$, for all $\Psi \subseteq \Phi$;\\[-.3in]
\item[c)] $G_{\Psi_1} \cap G_{\Psi_2} = G_{\Psi_1 \cup \Psi_2}$, for all $\Psi_1, \Psi_2 \subseteq \Phi$;\\[-.3in]
\item[d)] $\langle G_j \,|\, j \in I\rangle \cap \langle G_j \,|\, j\in J\rangle =
\langle G_j \,|\, j \in I\cap J\rangle$, for all $I, J \subseteq \{0, 1, 2\}$;\\[-.3in]
\item[e)] $G_0 \cdot G_2 = G_2 \cdot G_0 = G_{F_1}$. \\[-.3in]
\end{itemize}
\end{lemma}

\begin{proof}
These statements for stabilizer subgroups of $G$ are particular instances of similar such statements for flag-transitive subgroups of the automorphism group of a regular incidence complex of rank $3$ (or higher) obtained in \cite[\S2]{kom2}. The proof of the corresponding statements for the automorphism groups of polyhedra can also be found in \cite[pp.33,34]{rap}.

The proofs in \cite{kom2} rest on a crucial connectedness property (strong flag-connectedness) of incidence complexes, which in terms of $\K$ can be described as follows. Two flags of $\K$ are said to be {\em adjacent\/} ($j$-{\em adjacent\/}) if they differ in a single face (just their $j$-face, respectively). Here the number of flags $j$-adjacent to a given flag is $1$ if $j=1$ or $2$, or $r-1$ if $j=2$. Now the connectedness assumptions on $\K$ in Definition~\ref{compdef}(a,b) are equivalent to $\K$ being {\em strongly flag-connected\/}, in the sense that, if $\Upsilon$ and $\Omega$ are two flags of $\K$, then they can be joined by a sequence of successively adjacent flags $\Upsilon = \Upsilon_0,\Upsilon_1,\ldots,\Upsilon_k = \Omega$, each containing $\Upsilon\cap \Omega$. Thus, if $\Psi$ is a subset of the base flag $\Phi$, then any flag $\Omega$ containing $\Psi$ can be joined to $\Phi$ by a similar such sequence in which all flags contain $\Psi$. This key fact is the basis for an inductive argument. In fact, it can be shown by induction on $k$ that $\Omega$ is the image of $\Phi$ under an element of the subgroup $G_{\Psi}$ of $G$. Since the flag stabilizer $G_\Phi$ is a subgroup of $G_{\Psi}$, this then proves part $(b)$ and, in turn, parts (a), (c) and (d). Finally, part (e) reflects the fact that, for flags containing the base edge $F_1$, the two operations of taking $0$-adjacents or $2$-adjacents ``commute".
\end{proof}

Call an affine plane in $\E$ a {\em face mirror\/} of $\K$ if it contains a face of $\K$ and is the mirror of a reflection in $G(\K)$. Regular complexes with face mirrors must have planar faces. Moreover, if any one face of a regular complex $\K$ determines a face mirror, so do all faces, by the face-transitivity of $G(\K)$.

\begin{lemma}\label{reflectionface}
Let $\K$ be a regular complex with base flag $\Phi$ and group $G=G(\K)$. \\[-.3in]
\begin{itemize}
\item[a)] Then $G_\Phi$ has order at most $2$.\\[-.3in]
\item[b)] If $G_\Phi$ is non-trivial, then $\K$ has planar faces and the one non-trivial
element of $G_\Phi$ is the reflection in the plane containing the base face $F_2$; in particular, $\K$ has face mirrors.\\[-.3in]
\end{itemize}
\end{lemma}

\begin{proof}
Let $R \in G_\Phi$, that is, $F_i R = F_i$ for $i = 0, 1, 2$. Then $R$ fixes $F_0$, the
midpoint of $F_1$, and either the center of $F_2$ if $F_2$ is a finite polygon, or the axis of
$F_2$ if $F_2$ is an infinite polygon. Recall here that the faces of $\K$ are regular polygons. (The axis of a zigzag polygon is the line through the midpoints of its edges, and the axis of a helical polygon is the ``central" line of the ``spiral" staircase it forms.)

If the faces of $\K$ are finite, then $R$ fixes three non-collinear points, and hence $R$ must
be either trivial or the reflection in the plane spanned by the three points. In particular, the latter forces the faces to be planar.

If the faces are helical, then necessarily $R$ is trivial. In fact, since $R$ keeps the line through $F_1$ pointwise fixed and leaves the axis of $F_2$ invariant, this axis must necessarily intersect this line if $R$ is non-trivial. However, the latter cannot occur.

Finally, if the faces are zigzags, then $R$ again fixes the line through $F_1$ pointwise and leaves the axis of $F_2$ invariant. Since these two lines are not orthogonal to each other, $R$ must necessarily fix the plane containing $F_2$ pointwise, and, if not trivial, coincide with the reflection in this plane.

In either case, $G_\Phi$ has at most two elements. In particular, $G_\Phi$ can only be non-trivial if $F_2$ is planar and the plane through $F_2$ is a face mirror.
\end{proof}

In Section~\ref{facemir} we will explicitly characterize the regular complexes with a non-trivial flag stabilizer $G_\Phi$. Here we determine the structure of their subgroups $G_0, G_1$ and $G_2$.

\begin{lemma}\label{reflectiong2}
Let $\K$ be a regular complex with a non-trivial flag stabilizer $G_\Phi$ (of order $2$). Then,\\[-.3in]
\begin{itemize}
\item[a)] $G_0 \cong C_2 \times C_2 \cong G_1$; \\[-.3in]
\item[b)] $G_2\cong D_r$, the dihedral group of order $2r$. \\[-.3in]
\end{itemize}
\end{lemma}

\begin{proof}
Recall that the non-trivial element $R$ in $G_\Phi$ is the reflection in the plane containing the base face $F_2$. First note that the index of $G_\Phi$ in $G_0$ is $2$ (any element in $G_0$ either fixes $\Phi$ or moves $\Phi$ to the flag that is $0$-adjacent to $\Phi$). Suppose that $S \in G_0$ and $F_0 S \ne F_0$. Then $F_0 S^2 = F_0$ and hence $S^2 \in G_\Phi$. On the other hand, $S^2$ is a proper (orientation preserving) isometry and hence cannot coincide with $R$. Thus $S^2 = I$, the identity mapping. It follows that every element of $G_0$ has period $2$, so $G_0$ must be isomorphic to $C_2 \times C_2$. Similarly we also obtain $G_1 \cong C_2 \times C_2$, using $F_1$ instead of $F_0$.

Moreover, since $G_2$ fixes the line through $F_1$ pointwise and contains a reflection, namely $R$, the group $G_2$ must be dihedral of order $2r$. Note here that $G_2$ is transitive on the flags that are $2$-adjacent to $\Phi$, and that $R$ fixes $F_2$.
\end{proof}

The vertex-figures of a regular complex $\K$ with non-trivial flag-stabilizer are simple graphs on which the corresponding vertex-stabilizer subgroup of $G(\K)$ acts flag-transitively but not simply flag-transitively (the reflection in a face mirror of $\K$ acts trivially on a flag of the vertex-figure at a vertex of this face). Recall here that a flag of a graph (possibly with multiple edges) is a $2$-element set consisting of a vertex and an edge incident with this vertex.

We call a regular complex $\K$ {\em simply flag-transitive\/} if its {\em full\/} symmetry group $G(\K)$ is simply transitive on the flags of $\K$. Note here that our terminology strictly refers to a condition on the full group $G(\K)$, not a subgroup. Nevertheless, regular complexes $\K$ that are not simply flag-transitive in this sense, frequently admit proper subgroups of $G(\K)$ that do act simply flag-transitively on $\K$.

\begin{lemma}\label{r0r1}
Let $\K$ be a simply flag-transitive regular complex. Then,\\[-.3in]
\begin{itemize}
\item[a)] $G_0 = \langle R_0 \rangle$ and $G_1 = \langle R_1 \rangle$, for some point, line or plane reflection $R_0$ and some line or plane reflection $R_1$;\\[-.3in]
\item[b)] $G_2$ is a cyclic or dihedral group of order $r$ (in particular, $r$ is even if $G_2$ is dihedral);\\[-.3in]
\item[c)] $G_{F_2}=\langle R_0, R_1 \rangle$, a group isomorphic to a dihedral group (of finite or infinite order) and acting simply transitively on the flags of $\K$ containing $F_2$;\\[-.3in]
\item[d)] $G_{F_0}=\langle R_1, G_2 \rangle$, a finite group acting simply transitively on the flags of $\K$ containing~$F_0$.
\end{itemize}
\end{lemma}

\begin{proof}
Since $G_\Phi$ is trivial, the subgroups $G_0$ and $G_1$ are of order $2$ (recall here that $G_\Phi$ has index $2$ in $G_0$ and $G_1$) and hence are generated by involutions $R_0$ and $R_1$, respectively; the latter must necessarily be reflections in points, lines or planes. Note that $R_1$ cannot be a point reflection (in the point $F_0$), since otherwise this would force $F_2$ to be a linear apeirogon. Now parts (c) and (d) follow directly from Lemma~\ref{groupstruct}(b) (with $\Psi = \{F_2\}$ or $\Psi =\{F_0\}$, respectively). Here $G_{F_0}$ is finite, since $\K$ is discrete.

For part (b) note that $G_2$ keeps the line through $F_1$ pointwise fixed and hence is necessarily cyclic or dihedral. Since $G_2$ acts simply transitively on the flags of $\K$ containing $F_0$ and $F_1$, it must have order $r$.
\end{proof}

For a simply flag-transitive regular complex $\K$, the vertex-stabilizer subgroup $G_{F_0}$ of the base vertex $F_0$ in $G$ acts simply flag-transitively on (the graph that is) the vertex-figure of $\K$ at $F_0$. Its order is twice the number of faces of $\K$ containing $F_0$ (that is, twice the number of edges of the vertex-figure). We call $G_{F_0}$ the {\em vertex-figure group\/} of $\K$ at $F_0$. Even if the vertex-figure is a simple (geometric) graph, this group may {\it a priori\/} not be its full symmetry group in $\E$. If the vertex-figure has double edges, then any two edges with the same pair of end vertices can be mapped onto each other by an element of $G_{F_0}$ (swapping the two faces of $\K$ that determine these edges). On the other hand, the face stabilizer $G_{F_2}$ always turns out to be the full symmetry group of the face (face mirrors do not occur if $\K$ is simply flag-transitive).

The following theorem settles the enumeration of finite complexes. Recall that there are precisely eighteen finite regular polyhedra in $\E$, namely the Platonic solids and Kepler-Poinsot polyhedra, and their Petrie-duals (see \cite[Ch.7E]{rap}).

\begin{theorem}\label{fincomp}
The only finite regular complexes in $\E$ are the finite regular polyhedra.
\end{theorem}

\begin{proof}
First note that every finite subgroup of ${\mathcal O}(3)$ (the orthogonal group of $\E$) leaves a point invariant. Thus the center (centroid of the vertex set) of a finite regular complex $\K$ is invariant under $G$. Now assume to the contrary that $r \geq 3$. Then $G_2$ must contain a rotation about the line through $F_1$, which, in turn, forces the center of $\K$ to lie on this line. But then all lines through edges of $\K$ must pass through the center of $\K$, which is impossible.
\end{proof}

Thus we may restrict ourselves to the enumeration of infinite complexes. For the remainder of this section we will assume that $\K$ is a simply flag-transitive regular complex and that $\K$ is infinite. We continue to explore the structure of $G$.

First we say more about the stabilizer $G_{F_1}$ of the base edge $F_1$. In particular, $G_2$ has index $2$ in $G_{F_1}$, and
\begin{equation}\label{commut}
G_{F_1} \,=\, G_0 G_2 \,\cong\, G_2 \rtimes G_0.
\end{equation}
Note here that $G_0 \cap G_2$ is trivial, since $\K$ is simply flag-transitive. The semi-direct product in (\ref{commut}) is direct if $R_0$ is a point reflection (in the midpoint of $F_1$) or a plane reflection (in the perpendicular bisector of $F_1$).

The group $G$ of $\K$ is an infinite discrete group of isometries of $\E$. We say that
such a group acts {\em (affinely) irreducible} if there is no non-trivial linear subspace $L$ of
$\E$ which is invariant in the sense that $G$ permutes the translates of $L$ (and hence also the translates of the
orthogonal complement $L^\perp$). Otherwise, the group is called {\em (affinely) reducible}.

Our next theorem allows us to concentrate on regular complexes with an irreducible group. Note that $\K$ need not be simply flag-transitive for the purpose of this theorem.

\begin{theorem}\label{reducib}
A regular complex in $\E$ with an affinely reducible symmetry group is necessarily a regular polyhedron and hence is planar or blended in the sense of \cite[pp.221,222]{rap}.
\end{theorem}

\begin{proof}
We need to rule out the possibility that $r\geq 3$. Suppose to the contrary that $r \geq 3$ and that $G$ is reducible. Let $L$ be a plane in $\E$ whose translates are permuted under $G$ (this is the only case we need to consider). In particular, the images of $L$ under $G_2$ are translates of $L$. Since $G_2$ is cyclic or dihedral of order at least $3$, we know that its rotation subgroup is generated by a rotation $S$ about the line through $F_1$.

If $S$ is of order at least $3$, then its axis must necessarily be orthogonal to $L$. By taking a plane parallel to $L$ if need be, we may assume that $L$ is the perpendicular bisector of $F_1$. Then, by the edge-transitivity of $G$, each perpendicular bisector of an edge of $\K$ is an image of $L$ under $G$ and hence must be parallel to $L$. However, this forces all edges to be parallel, which is impossible (a complex cannot have linear apeirogons as faces).

Now let $S$ be of order $2$. The case when the axis of $S$ is orthogonal to $L$ can be
eliminated as before. It remains to consider the case when the axis of $S$ is parallel to $L$
and hence, without loss of generality, is contained in $L$.  Then, by the edge-transitivity of $G$, every edge of $\K$ is contained in an image of $L$ under $G$, and by the connectedness of $\K$, that image, being parallel to $L$, must necessarily be $L$ itself. Hence $\K$ must lie in a plane. Finally, to exclude this possibility note that a line segment in a plane can be an edge of at most two regular (finite or infinite) polygons of a given shape (consider the interior angles at vertices of the polygon). This forces $r=2$, contrary to our assumption.
\end{proof}

It is well known that an irreducible infinite discrete group of isometries in $\E$ is a crystallographic
group (that is, it admits a compact fundamental domain). The well-known Bieberbach theorem now tells us that such a group $G$ contains a subgroup $T(G)$ of the group ${\mathcal T}(\E)$ of translations of $\E$, such that the quotient $G/T(G)$ is finite; here we may think of $T(G)$ as a lattice in $\E$ (see \cite{bieber} and \cite[\S7.4]{ratcliffe}). If $R: x \mapsto xR' + t$ is a general element of $G$, with $R' \in {\mathcal O}(3)$ and $t \in \E$ a translation vector (we may thus think of $t \in {\cal T}(\E)$), then the mappings $R'$ clearly form a subgroup $G_*$ of ${\mathcal O}(3)$, called the {\em special group} of $G$. Thus $G_*$ is the image of $G$ under the epimorphism ${\mathcal I}(3)\mapsto {\mathcal O}(3)$ whose kernel is ${\cal T}(\E)$. (Here ${\mathcal I}(3)$ is the group of all isometries of $\E$.) In other words,
\[G_* \,=\, G\,{\cal T}(\E)/{\cal T}(\E) \,\cong\, G/(G \cap {\cal T}(\E)) \,=\, G/T(G),\]
if $T(G)$ is the full translation subgroup of $G$. In particular, $G_*$ is a finite group.

In essence, the following lemma was proved in \cite[p.220]{rap} (see also \cite[Lemma 3.1]{chiral2}). For the exclusion of rotations of period $6$ for groups in $\E$ observe that there is no finite subgroup of ${\mathcal O}(3)$ with two distinct $6$-fold axes of rotation, and that hence the presence of $6$-fold axes forces reducibility.

\begin{lemma}\label{casesnum}
The special group of an irreducible infinite discrete group of isometries in $\mathbb{E}^2$
or $\E$ does not contain rotations of periods other than $2, 3, 4$ or $6$. Moreover, a rotation of period $6$ can only occur in the planar case.
\end{lemma}

Lemma~\ref{casesnum} immediately limits the possible values of $r$. In fact, we must have $r = 2, 3, 4, 6$ or $8$. Similarly, the faces, if finite, must be $p$-gons with $p = 3, 4, 6$ or $8$. We later eliminate some of these possibilities.

The translation subgroup of the symmetry group (or any flag-transitive subgroup) of a non-planar regular complex in $\E$ is a $3$-dimensional lattice. The following lattices are particularly relevant for us. Let $a$ be a positive real number, let $k=1$, $2$ or $3$, and let ${\bf a} : = (a^k,0^{3-k})$, the vector with $k$ components $a$ and $3-k$ components $0$. Following \cite[p.166]{rap}, we write $\Lambda_{\bf a}$ for the sublattice of $a\mathbb{Z}^3$ generated by $\bf a$ and its images under permutation and changes of sign of coordinates. Observe that
\[  \Lambda_{\bf a} = a\Lambda_{(1^{k},0^{3-k})},  \]
when ${\bf a} = (a^{k},0^{3-k})$. Then $\Lambda_{(1,0,0)}=\mathbb{Z}^{3}$ is the standard {\em cubic lattice\/}; $\Lambda_{(1,1,0)}$ is the {\em face-centered cubic lattice\/} (with basis $(1,1,0)$, $(-1,1,0)$, $(0,-1,1)$) consisting of all integral vectors with even coordinate sum; and $\Lambda_{(1,1,1)}$ is the {\em body-centered cubic\/} lattice (with basis $(2,0,0)$, $(0,2,0)$, $(1,1,1)$).  

Occasionally we will appeal to the enumeration of the finite subgroups of ${\mathcal O}(3)$ (see \cite{grove}). In particular, under our assumptions on the complex $\K$, the special group $G_*$ must be one of five possible groups, namely
\begin{equation}
\label{finthree}
[3, 3], \;\; [3, 3]^+, \;\; [3, 3]^*, \;\; [3, 4] \;\; \mbox{or} \;\; [3, 4]^+.
\end{equation}
Here, $[p, q]$ and $[p, q]^+$, respectively, denote the full symmetry group or rotation subgroup of a Platonic solid $\{p, q\}$, and $[3, 3]^*$ is defined as
\[ [3, 3]^* \,=\, [3, 3]^+ \cup (-I)[3, 3]^+ \,\cong\, [3,3]^{+}\times C_2, \]
where the tetrahedron $\{3, 3\}$ is taken to be centered at the origin and $-I$ denotes the central inversion (generating $C_2$).
Note that $[3, 3]^*$ is the symmetry group of a common crystal of pyrite (known as fool's gold, due to its resemblance to gold).

We frequently use the following basic fact, which we record without proof.

\begin{lemma}\label{casesnum1}
Let $C$ be a cube, and let $R$ be the reflection in a plane determined by a pair of opposite face diagonals of $C$. If $S$ is any rotational symmetry of $C$ of period $4$ or $3$, respectively, whose rotation axis is not contained in the mirror of $R$, then $R$ and $S$ generate the full octahedral group $G(C)=[3,4]$ or its full tetrahedral subgroup $[3,3]$.
\end{lemma} 

Concluding this section, we also observe that the convex hull of the vertices adjacent to a given vertex $v$ of $\K$ is a vertex-transitive convex polygon or polyhedron, $P_v$ (say), in $\E$ and hence has all its vertices on a sphere centered at $v$. Here, $P_v$ is a polygon if and only if the vertex-figure of $\K$ at $v$ is planar; this can occur even if $\K$ is not a polyhedron (the complex $\K_{5}(1,2)$ of Section~\ref{com34} is an example). If $v = F_0$, then $G_2$ stabilizes the vertex of $P_v$ determined by the vertex of $F_1$ distinct from $F_0$.

\section{Complexes with face mirrors}
\label{facemir}

In this section we characterize the regular complexes with a non-trivial flag stabilizer $G_\Phi$ as 
the $2$-skeletons of regular $4$-apeirotopes in $\E$. We refer to \cite[Ch.7F]{rap} for the complete enumeration of these apeirotopes. Recall that there are eight regular $4$-apeirotopes in $\E$, two with finite $2$-faces and six with infinite $2$-faces (see (\ref{4apeirotopes}) below); they are precisely the discrete faithful realizations of abstract regular polytopes of rank $4$ in $\E$. 

Let $\K$ be a regular complex with group $G=G(\K)$, and let $G_\Phi$ be non-trivial, where again $\Phi$ denotes the base flag of $\K$. Then, by Theorem~\ref{reducib}, $G$ is (affinely) irreducible, since polyhedra do not have a non-trivial flag stabilizer. More importantly, by Lemma~\ref{reflectionface}, $\K$ has face mirrors. In particular, the faces of $\K$ are planar and the only non-trivial element of $G_\Phi$ is the reflection $R_3$ (say) in the plane containing $F_2$. Moreover, by Lemma~\ref{reflectiong2}, $G_2$ is a dihedral group of order $2r$ containing $R_3$.

Recall from Lemma~\ref{reflectiong2}(a) that $G_0$ and $G_1$ each contain two involutions distinct from $R_3$ whose product is $R_3$. Hence, since $R_3$ is an improper isometry, one must be a proper isometry (a half-turn) and the other an improper isometry (a point reflection or plane reflection). Let $R_0$ and $R_1$, respectively, denote the improper isometries in $G_0$ and $G_1$ distinct from $R_3$.

First we discuss complexes $\K$ with finite faces. Since $G_0$ fixes the centers of $F_1$ and $F_2$, the isometry $R_0$ cannot be a point reflection and hence must be the plane reflection in the perpendicular bisector of $F_1$. Similarly, since $G_1$ fixes $F_0$ and the center of $F_2$, the isometry $R_1$ is also a plane reflection.

Let $R_2$ be a plane reflection such that $R_2, R_3$ are distinguished generators for the dihedral group $G_2 \cong D_r$ with mirrors inclined at $\pi/r$. Then
\[G = \langle G_0, G_1, G_2 \rangle = \langle R_0, R_1, R_2, R_3 \rangle,\]
and $G$ is a discrete irreducible reflection group in $\E$. We know that the generator $R_3$ commutes with $R_0$ and $R_1$ by Lemma~\ref{reflectiong2} (a). Moreover, the mirrors of $R_0$ and $R_2$ are perpendicular, and hence $R_0$ and $R_2$ commute as well. It follows that the Coxeter diagram associated with the generators $R_0, R_1, R_2, R_3$ of $G$ is a string diagram. In particular, $G$ must be an infinite group (geometrically) isomorphic to the Coxeter group $[4, 3, 4]$, the symmetry group of the regular cubical tesellation $\{4,3,4\}$ in $\E$ (see \cite{coxeter}). Note here that the subgroups $\langle R_0, R_1, R_2 \rangle$ and $\langle R_1, R_2, R_3 \rangle$ each are finite crystallographic (plane) reflection groups in $\E$.

We now employ a variant of Wythoff's construction to show that $\K$ coincides with the $2$-skeleton of the cubical tessellation $\{4,3,4\}$ of $\E$ associated with the generators $R_0, R_1, R_2, R_3$.

First note that the base vertex $F_0$ is fixed by $G_1$ and $G_2$, and hence by $R_1, R_2$ and $R_3$. This implies that $F_0$ is just the base (initial) vertex of $\{4, 3, 4\}$ and that the subgroup generated by $R_1, R_2, R_3$ is its stabilizer in $G$. In $\K$, the orbit of $F_0$ under $G_0 = \langle R_0, R_3 \rangle$ consists only of $F_0$ and $F_0 R_0$, so the base edge $F_1$ of $\K$ is the line segment joining $F_0$ and $F_0 R_0$. Hence $F_1$ is just the base edge of $\{4, 3, 4\}$. Similarly, the orbit of $F_1$ under $\langle G_0, G_1 \rangle = \langle R_0, R_1, R_3 \rangle$ in $\K$ coincides with the orbit of $F_1$ under the subgroup $\langle R_0, R_1 \rangle$, and hence the base face $F_2$ of $\K$ is just the base face of $\{4, 3, 4\}$. Finally, then, the vertices, edges and faces of $\K$ are just the images of $F_0$, $F_1$ and $F_2$ under $G$, and hence $\K$ is the $2$-skeleton of $\{4, 3, 4\}$.

Before proceeding with complexes $\K$ with infinite faces, we briefly review the {\em free abelian apeirotope\/} or ``apeir" construction described in \cite{pm,pm1}. This construction actually applies to any finite (rational) regular polytope $\mathcal Q$ in some euclidean space, but here we only require it for regular polyhedra $\mathcal Q$ in $\E$, where it produces a regular apeirotope of rank $4$ in $\E$ with vertex-figure $\mathcal Q$.  

Let $\mathcal Q$ be a finite regular polyhedron in $\E$ with symmetry group 
$G({\mathcal Q}) = \langle \widehat{R}_1,\widehat{R}_2,\widehat{R}_3\rangle$ (say), where the labeling of the distinguished generators begins at $1$ deliberately. Let $o$ be the centroid of the vertex-set of $\mathcal Q$, let $w$ be the initial vertex of $\mathcal Q$, and let $\widehat{R}_0$ denote the reflection in the point $\frac{1}{2}w$. Then there is a regular $4$-apeirotope in $\E$, denoted $\apeir {\mathcal Q}$, with $\widehat{R}_0,\widehat{R}_1,\widehat{R}_2,\widehat{R}_3$ as the generating reflections of its symmetry group, $o$ as initial vertex, and $\mathcal Q$ as vertex-figure. In particular, $\apeir Q$ is discrete if $\mathcal Q$ is rational (the vertices of $\mathcal Q$ have rational coordinates with respect to some coordinate system). The latter limits the choices of $\mathcal Q$ to $\{3,3\}$, $\{3,4\}$ or $\{4,3\}$, or their Petrie duals $\{4,3\}_3$, $\{6,4\}_3$ or $\{6,3\}_4$, respectively; the six corresponding $4$-apeirotopes are listed in (\ref{4apeirotopes}).

Now consider complexes $\K$ with infinite (planar) faces, that is, zigzag faces. Now the generator $R_0$ must be a point reflection in the midpoint of $F_1$. In fact, if $R_0$ was a plane reflection (necessarily in the perpendicular bisector of $F_1$), the invariance of $F_1$ and $F_2$ under $R_0$ would force $F_2$ to be a linear apeirogon, which is not possible. On the other hand, $R_1, R_2$ and $R_3$ still are plane reflections, with $R_2$ and $R_3$ as before. In fact, if $R_1$ was a point reflection (necessarily in the point $F_0$), then again $F_2$ would have to be a linear apeirogon. Moreover, since $R_1$ leaves $F_0$ and $F_2$ invariant, its plane mirror must necessarily be perpendicular to the plane containing $F_2$. Hence, since the latter is the mirror of $R_3$, the generators $R_1$ and $R_3$ must commute (this also follows from Lemma~\ref{reflectiong2}(a)). Similarly, the point reflection $R_0$ commutes with $R_2$ and $R_3$, because its fixed point (the midpoint of $F_1$) lies on the mirrors of $R_2$ and $R_3$. In addition, the generators $R_1$, $R_2$, $R_3$ all fix $F_0$, so that $\langle R_1, R_2, R_3 \rangle$ is a finite, irreducible crystallographic reflection group in $\E$. This group has a string Coxeter diagram;
note here that $R_1$ and $R_2$ cannot commute, again since faces are not linear apeirogons. Hence, this group must necessarily be $[3, 3]$, $[3, 4]$ or $[4, 3]$, with $R_1, R_2, R_3$ the distinguished generators.

Now we can prove that $\K$ is the $2$-skeleton of the regular $4$-apeirotope $\apeir\{3,3\}$, $\apeir\{3,4\}$ or $\apeir\{4,3\}$ in $\E$. For convenience we take $F_0$ to be the origin $o$, so that the subgroup $\langle R_1, R_2, R_3 \rangle$ is just the standard Coxeter group $[3, 3], [3, 4]$ or $[4, 3]$ with distinguished generators $R_1, R_2, R_3$. Set $F_0':=F_0 R_0$, so that $F_1$ has vertices $F_0$ and $F_0'$. Then $R_0$ is the reflection in the point $\frac{1}{2} F_0'$, and the orbit of $F_0'$ under $\langle R_1, R_2, R_3 \rangle$ is just the vertex set of a copy of ${\mathcal Q}=\{3, 3\}$, $\{3, 4\}$ or $\{4, 3\}$, respectively. Thus the configuration of the mirrors of $R_0, R_1, R_2, R_3$ is exactly that of the mirrors of the distinguished generators for the symmetry group of the $4$-apeirotope $\apeir{\mathcal Q}$. Hence $G$ must be the symmetry group of $\apeir{\mathcal Q}$. Moreover, the initial vertex, $F_0 = o$, is the same in both cases, so the vertex sets of $\K$ and $\apeir{\mathcal Q}$ must be the same as well. The proof that the edges and faces of $\K$ are just those of $\apeir{\mathcal Q}$ follows from the same general argument as for complexes with finite faces.

In summary we have established the following theorem.

\begin{theorem}
\label{nontriv}
Every regular polygonal complex with a non-trivial flag stabilizer is the
$2$-skeleton of a regular $4$-apeirotope in $\E$.
\end{theorem}

Let $\mathcal P$ be a regular $4$-apeirotope in $\E$, and let $\pi$ denote the Petrie operation on (the vertex-figure of) $\mathcal P$. Recall that, if $\Gamma({\mathcal P}) = \langle T_0, T_1, T_2, T_3 \rangle$, then $\pi$ is determined by changing the generators on $\Gamma({\mathcal P})$ according to
\[\pi : \, (T_0, T_1, T_2, T_3) \mapsto (T_0, T_1T_3, T_2, T_3).\]
The $4$-apeirotope associated with the new generators on the right is called the Petrie-dual of
$\mathcal P$ and is denoted by ${\mathcal P}^\pi$ (see \cite[Ch.7F]{rap}).

Our above analysis shows that $\K$ is the $2$-skeleton of either the cubical tessellation
$\{4, 3, 4\}$ if the faces are finite, or the regular $4$-apeirotope $\apeir\{3,3\}$, $\apeir\{3,4\}$ or $\apeir\{4,3\}$ if the faces are infinite. This covers four of the eight regular $4$-apeirotopes in $\E$. However, the eight apeirotopes occur in pairs of Petrie-duals and, as we show in the next lemma, the Petrie operation does not affect the $2$-skeleton. Thus, all eight apeirotopes actually occur but contribute only four regular complexes.

\begin{lemma}\label{petrieop}
If $\mathcal P$ is a regular $4$-apeirotope in $\E$, then $\mathcal P$ and ${\mathcal P}^\pi$ have the
same $2$-skeleton.
\end{lemma}

\begin{proof}
First recall that, for regular $3$-polytopes, the Petrie-dual has the same vertices and edges as the original polytope. Since the vertex-figure of ${\mathcal P}^\pi$ is just the Petrie-dual of the vertex-figure of $\mathcal P$, the base vertices of $\mathcal P$ and ${\mathcal P}^\pi$ are the same (they are the fixed points of the vertex-figure groups), and hence the vertex sets of $\mathcal P$ and ${\mathcal P}^\pi$ are the same. In a $4$-apeirotope, the edges and $2$-faces containing a given vertex correspond to the vertices and edges of the vertex-figure at that vertex. Therefore, since the vertex-figures of $\mathcal P$ and ${\mathcal P}^\pi$ at their base vertices have the same vertices and edges, the sets of edges and $2$-faces of $\mathcal P$ and ${\mathcal P}^\pi$ containing the base vertex must necessarily also be the same. It follows that $\mathcal P$ and ${\mathcal P}^\pi$ have the same $2$-skeleton.
\end{proof}

Concluding this section we list the eight regular $4$-apeirotopes in a more descriptive way in pairs of Petrie-duals using the notation in \cite{rap}. The apeirotopes in the top row have square faces, and their facets are cubes or Petrie-Coxeter polyhedra $\{4, 6 \,|\,4\}$. All others have zigzag faces, and their facets are blends of the Petrie-duals $\{\infty,3\}_6$ or $\{\infty,4\}_4$ of the plane tessellations $\{6,3\}$ or $\{4,4\}$, respectively, with the line segment $\{ \,\}$ or linear apeirogon $\{\infty\}$ (see \cite[Ch.7E]{rap}).
\begin{equation}\label{4apeirotopes}
\begin{array}{ccc}
\{4, 3, 4\} & &\{\{4, 6 \,|\,4\}, \{6, 4\}_3\}\\[.04in]
\apeir \{3, 3\} = \{\{\infty, 3\}_6 \# \{ \, \}, \{3, 3\}\} 
& &\{\{\infty, 4\}_4 \# \{\infty\}, \{4, 3\}_3\} = \apeir\{4, 3\}_3 \\[.04in]
\apeir \{3, 4\} = \{\{\infty, 3\}_6 \# \{ \, \}, \{3, 4\}\} 
& &\{\{\infty, 6\}_3 \# \{\infty\}, \{6, 4\}_3\} = \apeir\{6, 4\}_3 \\[.04in]
\apeir \{4, 3\} = \{\{\infty, 4\}_4 \# \{\, \}, \{4, 3\}\} 
& &\{\{\infty, 6\}_3 \# \{\infty\}, \{6, 3\}_4\} = \apeir \{6, 3\}_4\\[.04in]
\end{array}
\end{equation}
Note here that $(\apeir {\mathcal Q})^\pi = \apeir({\mathcal Q}^\pi)$ for each vertex-figure $\mathcal Q$ that occurs. The parameter $r$ counting the number of faces around an edge of the $2$-skeleton $\K$ is just the last entry in the Schl\"afli symbol (the basic symbol $\{p,q,r\}$) of the corresponding $4$-apeirotope (or, equivalently, its Petrie dual). Hence, $r=4$, $3$, $4$ or $3$, respectively.

In summary we have the following theorem.

\begin{theorem}
\label{Petskels}
Up to similarity, there are precisely four regular polygonal complexes with a non-trivial flag stabilizer, each given by the common $2$-skeleton of two regular $4$-apeirotopes in $\E$ which are Petrie duals of each other.
\end{theorem}

Now that the regular complexes with non-trivial flag stabilizers have been described, we can restrict ourselves to the enumeration of simply flag-transitive regular complexes.

\section{Operations}
\label{opera}

In this section we discuss operations on the generators of the symmetry group of a regular polygonal complex that allow us to construct new complexes from old. In particular, this will help reduce the number of cases to be considered in the classification.

Let $\K$ be a regular complex in $\E$, and let $\K$ have a simply flag-transitive (full symmetry) group $G = G(\K) = \langle G_0, G_1, G_2 \rangle$. Then recall from Lemma~\ref{r0r1}(a) that 
$G_0 = \langle R_0 \rangle$ for some point, line or plane reflection $R_0$, that $G_1 = \langle R_1 \rangle$ for some line or plane reflection $R_1$, and that $G_2$ is cyclic or dihedral of order $r$.

We will define two operations on the group $G$ which replace $R_0$ or $R_1$, respectively,
but retain $G_2$. They employ elements $R$ of $G_{2}$ with the property that $R_0 R$ or $R_1 R$, respectively, is an involution. Here we are not assuming that $R$ itself is an involution; however, this will typically be the case.

\subsection{Operation $\lambda_0$}
\label{opl0}

To begin with, suppose $R$ is an element of $G_2$ such that $R_0 R$ is an involution (if $R_0$ is a point or plane reflection, this actually forces $R$ to be an involution as well). At the group level we define our first operation $\lambda_0$ on the generating subgroups $G_{0}=G_{0}(\K)$, $G_{1}=G_{1}(\K)$, $G_{2}=G_{2}(\K)$ of $G$ by way of
\begin{equation}\label{opone}
\lambda_0 = \lambda_0(R)\!:\;\,  (R_0, R_1, G_2)\; \mapsto\; (R_0 R, R_1, G_2).
\end{equation}

\begin{lemma}\label{lambda1}
Let $\K$ be a simply flag-transitive regular complex with group $G = \langle R_0, R_1, G_2 \rangle$, and let $R$ be an element in $G_2$ such that $R_0 R$ is an involution. Then 
there exists a regular complex, denoted $\K^{\lambda_0}$, with the same vertex set and edge set as $\K$ and with the same symmetry group $G$, such that
\begin{equation}\label{klambda}
\langle R_0 R \rangle \subseteq G_0(\K^{\lambda_0}), \quad
  G_1(\K) = \langle R_1 \rangle \subseteq G_1(\K^{\lambda_0}),\quad
\,\, G_2(\K) = G_2(\K^{\lambda_0}).
\end{equation}
The inclusions in (\ref{klambda}) are equalities if and only if $\K^{\lambda_0}$ is simply flag-transitive.
\end{lemma}

\begin{proof}
Let $R_0' := R_0 R$, and let $G_0' := \langle R_0 R \rangle$, $G_1' := G_1$ and $G_2' := G_2$. We shall obtain the complex $\K^{\lambda_0}$ by Wythoff's construction (see \cite[Ch.5B]{rap}). First note that the subgroup $\langle G_1', G_2' \rangle$ has exactly one fixed point, namely the base vertex $F_0$ of $\K$. Thus we also take $F_0$ as base (initial) vertex for $\K^{\lambda_0}$. Since both $R_0$ and $R$ leave the base edge $F_1$ of $\K$ invariant, the orbit of $F_0$ under $G_0'$ consists of the two vertices of $F_1$, and hence $F_1$ can also serve as base edge of $\K^{\lambda_0}$. Moreover, the vertex set and edge set, respectively, of the base face $F_2'$ of $\K^{\lambda_0}$ are determined by the orbits of $F_0$ and $F_1$ under the subgroup $\langle R_0', R_1 \rangle$. Finally, then, the complex $\K^{\lambda_0}$ consists of all the vertices, edges and faces obtained as images of $F_0$, $F_1$ and $F_2'$ under $G$.

Note that $\K^{\lambda_0}$ really is a complex. By construction, the graphs of $\K$ and $\K^{\lambda_0}$ coincide, since their vertex sets and edge sets are the same. In particular, the graph of $\K^{\lambda_0}$ is connected. Moreover, (possibly) up to a change of edge multiplicities (from single to double, or vice versa), the vertex-figures of $\K$ and $\K^{\lambda_0}$ at the common base vertex $F_0$ coincide, since each can be obtained by Wythoff's construction applied to $\langle R_1, G_2 \rangle$ with base vertex $F_0 R_0 = F_0 R_0'$. In particular, the vertex-figures remain connected. The discreteness follows from our comments made after Definition \ref{defcomplex}. Finally, by the edge-transitivity of $G$ on $\K^{\lambda_0}$ and the fact that $G_2$ remained unchanged, each edge must again be contained in a fixed number of faces; this number is $r$ if $\K^{\lambda_0}$ is simply flag-transitive.
\end{proof}

Notice that the operation $\lambda_0$ actually also applies (with slight modifications) to regular complexes whose (full symmetry) group is not simply flag-transitive, that is, to the $2$-skeletons of regular $4$-apeirotopes. In this more general context, the operation $\lambda_0$ becomes invertible and its inverse is associated with the element $R^{-1}$ of $G_2$. There are instances in Lemma~\ref{lambda1} when the group of the resulting complex $\K^{\lambda_0}$ is no longer
simply flag-transitive, and then the more general operation is needed to recover the original complex from it.

For example, if $\K$ is the regular complex whose faces are all the Petrie polygons of all the cubes of the cubical tessellation $\{4, 3, 4\}$ of $\E$ (and the vertices and edges of $\K$ are just those of $\{4, 3, 4\}$), then for a suitable choice of $R$ in $G_2$, the new complex $\K^{\lambda_0}$ is the $2$-skeleton of $\{4, 3, 4\}$. In the notation of Section~\ref{com34}, $\K=\K_{6}(1,2)$. In this example, the vertex-figures of $\K^{\lambda_0}$ are simple graphs, while those of $\K$ have double edges.

\subsection{Operation $\lambda_1$}
\label{opl1}

Next we consider the operation $\lambda_1$ on the subgroups $G_0$, $G_1$ and $G_2$ of $G$ defined by
\begin{equation}\label{optwo}
\lambda_1 = \lambda_1(R) : \:\, (R_0, R_1, G_2)\; \mapsto\; (R_0, R_1R, G_2),
\end{equation}
where now $R \in G_2$ is such that $R_1 R$ is an involution. Again we shall concentrate on complexes with a simply flag-transitive group, although the operation applies (with slight modifications) more generally to arbitrary regular complexes. Now we have

\begin{lemma}\label{lambda2}
Let $\K$ be a simply flag-transitive regular complex with group $G = \langle R_0, R_1, G_2 \rangle$, and let $R$ be an element in $G_2$ such that $R_1 R$ is an involution. Then 
there exists a regular complex, denoted $\K^{\lambda_1}$ and again simply flag-transitive, 
with the same vertex set and edge set as $\K$ and with the same group $G$, such
that
\begin{equation}\label{klambda1}
G_0(\K) = \langle R_0 \rangle = G_0(\K^{\lambda_1}), \quad
  \langle R_1 R \rangle = G_1(\K^{\lambda_1}),\quad
\,\, G_2(\K) = G_2(\K^{\lambda_0}).
\end{equation}
\end{lemma}

\begin{proof}
The proof follows the same pattern as the proof of Lemma \ref{lambda1}. As before, we employ Wythoff's construction with initial vertex $F_0$ to construct $\K^{\lambda_1}$. Again, the edge graphs of $\K$ and $\K^{\lambda_1}$ are the same; note here that the element $R_1' := R_1 R = R^{-1} R_1$ maps $F_1$ to $F_1 R_1$, which is the other edge of $F_2$ containing $F_0$.

The connectedness and discreteness of $\K^{\lambda_1}$ are derived as before. The vertex-figure of $\K^{\lambda_1}$ at $F_0$ is obtained from Wythoff's construction with group $\langle R_1 R, G_2 \rangle = \langle R_1, G_2 \rangle$ and initial vertex $F_0 R_0$, and hence is the same as the vertex-figure of $\K$ at $F_0$, again (possibly) up to a change of edge multiplicities. In particular, the vertex-figures remain connected.

However, unlike in the previous case, the complex $\K^{\lambda_1}$ cannot have a non-trivial flag stabilizer. In fact, when applied to a regular complex with a non-trivial flag stabilizer $\langle R_3 \rangle$, the operation $\lambda_1$ must necessarily be the Petrie operation (on the vertex-figure) of the corresponding regular $4$-apeirotope (there is just one possible choice for $R$, namely $R_3$) and hence, by Lemma \ref{petrieop}, must leave the $2$-skeleton invariant. Now suppose $\K^{\lambda_1}$ has a non-trivial flag stabilizer and hence is the $2$-skeleton of a regular $4$-apeirotope. Then, since $\lambda_1$ is invertible, $\K$ can be recovered from $\K^{\lambda_1}$ by the inverse of $\lambda_1$, which necessarily must be the Petrie operation on (the apeirotope
determined by) $\K^{\lambda_1}$. It follows that the original complex $\K$ must have had a non-trivial flag stabilizer as well. Thus $\K^{\lambda_1}$ is simply flag-transitive.

Finally, by Lemma~\ref{r0r1}(b), the number of faces around an edge of $\K^{\lambda_1}$ is $r$, since $G_2(\K) = G_2(\K^{\lambda_1})$.
\end{proof}

\subsection{Mirror vectors $(2, k)$ from $(0, k)$ by way of $\lambda_0$}
\label{m2k0kl0}

We now discuss some particularly interesting applications of $\lambda_0$ and $\lambda_1$. Once again we restrict ourselves to regular complexes $\K$ with simply flag-transitive symmetry groups. We call the vector $(dim(R_0), dim(R_1))$ the {\em mirror vector} of $\K$, where $dim(R_i)$ is the dimension of the mirror (fixed point set) of the symmetry $R_i$ for $i=0,1$. Here, to keep notation simple, we suppress information about the subgroup $G_2$, although this is not to imply that the structure of $G_2$ is irrelevant. We sometimes indicate the mirror vector explicitly and denote a complex $\K$ with mirror vector $(j, k)$ by $\K(j, k)$.

If $\K$ is a regular polyhedron, then $G_2$ is generated by a (point, line or plane) reflection $R_2$ and we refer to the vector $(dim(R_0), dim(R_1),dim(R_2))$ as the {\em complete mirror vector} of $\K$. (In \cite[Ch.7E]{rap}, this vector was called the dimension vector.) In the  context of the present paper we usually take the enumeration of regular polyhedra for granted and concentrate on complexes that are not polyhedra. 

We begin with the operation $\lambda_0$. Let $\K$ be a regular complex with mirror vector $(2, k)$ for some $k = 1, 2$; then $R_0$ is the reflection in the perpendicular bisector of $F_1$. (Recall here that the case $k = 0$ cannot occur by Lemma~\ref{r0r1}(a).) If $R \in G_2$ is a half-turn, then its mirror (axis) must necessarily be the line through $F_1$, so $R_0 R$ must necessarily be the point reflection in the midpoint of $F_1$. It follows that the corresponding
complex $\K^{\lambda_0}$ has mirror vector $(0, k)$. Conversely, let $\K$ be a regular 
complex with mirror vector $(0, k)$ for some $k = 1, 2$. If $R \in G_2$ is a half-turn, then its mirror must necessarily contain the mirror of $R_0$ (the midpoint of $F_1$), so $R_0 R$ must necessarily be the plane reflection in the perpendicular bisector of $F_1$. Hence the corresponding complex $\K^{\lambda_0}$ has mirror vector $(2, k)$. In either case, $\lambda_0$ is an involutory operation.

The next lemma allows us to derive the enumeration of regular polygonal complexes with
$dim(R_0) = 2$ from the enumeration of those with $dim(R_0) = 0$.

\begin{lemma}\label{appone}
Let $\K$ be an infinite simply flag-transitive regular complex with an affinely irreducible group and mirror vector $(2, k)$ for some $k = 1, 2$. Then $G_2(\K)$ contains a half-turn $R$. In particular, the corresponding complex $\K^{\lambda_0}$, with $\lambda_0 = \lambda_0(R)$, is a regular complex with mirror vector $(0, k)$, and $\K = (\K^{\lambda_0})^{\lambda_0}$.
\end{lemma}

\begin{proof}
Appealing to the enumeration of finite subgroups of isometries of $\E$, we first observe that the special subgroup $G_*$ of $G$ must be a subgroup of the octahedral group $[3, 4]$ (see (\ref{finthree})); note here that $G_*$ is linearly irreducible and that $G$ is an infinite group. Now suppose that $G_2$ does not contain a half-turn. Then $G_2$ must be a cyclic group $C_3$ or dihedral group $D_3$ and hence must contain a $3$-fold rotation with axis perpendicular to the plane mirror of $R_0$. On the other hand, the octahedral group does not contain a $3$-fold rotation with an axis perpendicular to a plane reflection mirror. Thus $G_2$ must contain a half-turn. Now the lemma follows from our previous considerations.
\end{proof}

Recall from Theorem~\ref{reducib} that a regular complex in $\E$ with a reducible group must be a planar or blended polyhedron. Of these polyhedra, only the three planar tessellations $\{3,6\}$, $\{4,4\}$ and $\{6,3\}$ have mirror vectors of type $(2,k)$ (in fact, of type $(2,2)$), and then $G_2$ is generated by a plane reflection and does not contain a half-turn. In other words, our irreducibility assumption in Lemma~\ref{appone} really only eliminates these three choices for $\K$.
 
Furthermore, notice for Lemma~\ref{appone} that, vice versa, we may not generally be able to similarly obtain a given regular complex with mirror vector $(0, k)$ from one with vector $(2, k)$. In fact, the octahedral group does contain $3$-fold rotations about axes passing through the mirror of a point reflection (namely, the central involution), so $G_{2}$ might not contain a half-turn.

\subsection{Mirror vector $(0, k)$ from $(1, k)$ by way of $\lambda_0$}
\label{m0k1kl0}

Now let $\K$ be a regular complex with mirror vector $(0, k)$ for some $k = 1, 2$ (as before, the case $k=0$ cannot occur); then $R_0$ is the point reflection in the midpoint of $F_1$.  If $R \in G_2$ is a plane reflection, then its mirror must necessarily contain $F_1$ (and hence the mirror of $R_0$), so $R_0 R$ must be the half-turn about the line perpendicular to the mirror of $R$. Then the corresponding
complex $\K^{\lambda_0}$ has mirror vector $(1, k)$. Conversely, let $\K$ be a regular 
complex with mirror vector $(1, k)$ for some $k = 1, 2$; then $R_0$ is a half-turn about an axis perpendicular to the line through $F_1$ and passing through the midpoint of $F_1$. If $R \in G_2$ is a plane reflection whose mirror is perpendicular to the mirror of $R_0$, 
then $R_0 R$ is a point reflection in the midpoint of $F_1$. Now the corresponding complex $\K^{\lambda_0}$ has mirror vector $(0, k)$. In either case, $\lambda_0$ is involutory.

The following lemma allows us to deduce the enumeration of regular complexes with
$dim(R_0) = 0$ and a dihedral subgroup $G_2$ from the enumeration of those with $dim(R_0) = 1$ (and a dihedral group $G_2$).

\begin{lemma}\label{apptwo}
Let $\K$ be a simply flag-transitive regular complex with an affinely irreducible group, a dihedral subgroup $G_2$, and mirror vector $(0, k)$ for some $k = 1, 2$. Then, for any plane reflection $R \in G_2$, the corresponding complex $\K^{\lambda_0}$,
with $\lambda_0 = \lambda_0(R)$, is a regular complex with mirror vector $(1, k)$. In
particular, $\K = (\K^{\lambda_0})^{\lambda_0}$.
\end{lemma}

\begin{proof}
All we need to say here is that $R_0 R$ is the half-turn about the line through the midpoint of $F_1$ perpendicular to the mirror of $R$. Then Lemma~\ref{lambda1} applies.
\end{proof}

As in the previous case, we may not generally be able to obtain a given regular complex
with mirror vector $(1, k)$ from a complex with vector $(0, k)$, unless the axis of the half-turn $R_0$ is perpendicular to the mirror of a plane reflection in $G_2$. However, the latter condition is not guaranteed.

It is not hard to see that besides the scenarios described in Lemmas~\ref{appone} and \ref{apptwo} there are only three other possible choices of $R$ in which $\lambda_0$ can be applied. First, if $\K$ has mirror vector $(2, k)$ and $G_2$ is dihedral, then trivially any plane
reflection $R$ in $G_2$ has mirror perpendicular to the mirror of $R_0$ and leads to
a new complex $\K^{\lambda_0}$ with mirror vector $(1, k)$. Second, and conversely, if $\K$ has mirror vector $(1, k)$, $G_2$ is dihedral, and $G_2$ contains a plane reflection $R$ whose mirror contains the axis of $R_0$, then we arrive at a new complex $\K^{\lambda_0}$ with mirror vector $(2, k)$. For these two choices, $\K = (\K^{\lambda_0})^{\lambda_0}$. Third, if $\K$ has mirror vector $(1, k)$ and $R$ is a rotation in $G_2$, then $R_0 R$ is the half-turn about an axis contained in the perpendicular bisector of $F_1$, and hence $\K^{\lambda_0}$ is again a complex with mirror vector $(1, k)$. For the third choice, note that $R$ need not be an involution and thus $\lambda_0$ need not be involutory.

\subsection{Mirror vector $(k, 1)$ from $(k, 2)$ by way of $\lambda_1$}
\label{mk1k2l1}

Next we turn to applications of $\lambda_1$. Let $\K$ be a regular complex with affinely irreducible symmetry group $G$, and let $\K$ be simply flag-transitive. Assume that the mirror vector of $\K$ is of the form $(k, 1)$ for some $k = 0, 1, 2$, and that $G_2$ contains a plane reflection $R$ whose mirror contains the axis of the half-turn $R_1$. Then $R_1 R$ is the reflection in the plane that is perpendicular to the mirror of $R$ and contains the axis of $R_1$. It follows that the corresponding regular complex $\K^{\lambda_1}$ has mirror vector $(k, 2)$. Conversely, if $\K$ has mirror vector $(k, 2)$ and $G_2$ contains a plane reflection $R$ whose mirror is perpendicular to the plane mirror of $R_1$, then $R_1 R$ is the half-turn whose axis is the intersection of the two mirrors. Now the complex $\K^{\lambda_1}$ has mirror vector $(k, 1)$. Once again, in either case the operation $\lambda_1$ is involutory.

In summary, appealing to Lemma~\ref{lambda2}, we now have the following result that allows us to derive the enumeration of certain regular complexes with $dim(R_1) = 1$ from the enumeration of those with $dim(R_1) = 2$.

\begin{lemma}\label{appthree}
Let $\K$ be a simply flag-transitive regular complex with an affinely irreducible group
and mirror vector $(k, 1)$ for some $k = 0, 1, 2$. Assume also that $G_2$ contains a plane reflection $R$ whose mirror contains the axis of the half-turn $R_1$. Then the corresponding complex $\K^{\lambda_1}$, with $\lambda_1 = \lambda_1(R)$, is a regular complex with mirror vector $(k, 2)$. In particular, $\K = (\K^{\lambda_1})^{\lambda_1}$.
\end{lemma}

It can be shown that the above analysis exhausts all possible choices for elements $R$ in $G_2$  that can lead to a new complex $\K^{\lambda_1}$ via the corresponding operations $\lambda_{1}=\lambda_{1}(R)$.

\section{Complexes with finite faces and mirror vector $(1, 2)$}
\label{finmir12}

In this section and the next, we enumerate the simply flag-transitive regular complexes with mirror vector $(1, 2)$. Here we begin with complexes with finite faces. Throughout, let $\K$ be an infinite simply flag-transitive regular complex with an affinely irreducible group $G = \langle R_0, R_1, G_2 \rangle$.

\subsection{The special group}
\label{specgroup}

As we remarked earlier, if a regular complex $\K$ has finite faces, then these faces are necessarily planar or skew regular polygons. If additionally $\K$ has mirror vector $(1,2)$, then $R_0$ is a line reflection (half-turn) and $R_1$ a plane reflection whose mirrors intersect at the center of the base face $F_2$. Let $E_1$ denote the perpendicular bisector of the base edge $F_1$. Then the subgroup $G_2$ leaves the line through $F_1$ pointwise fixed and acts faithfully on $E_1$. The product $R_0 R_1$ is a rotatory reflection that leaves the plane $E_2$ through the midpoints of the edges of $F_2$ invariant, and its period is just the number of vertices of $F_2$. Note that $E_2$ is perpendicular to the plane mirror of $R_1$ and necessarily contains the axis of $R_0$. Since the axis of $R_0$ also lies in $E_1$ and $E_1 \ne E_2$, it must necessarily coincide with $E_1 \cap E_2$.

Throughout, it is convenient to distinguish the following two alternative scenarios for $R_0$ and $G_2$ and treat them as separate cases.\\[-.3in]
\begin{itemize}
\item [(A)] The axis of $R_0$ is contained in the mirror of a plane reflection, $R_2$ (say), in $G_2$.\\[-.3in]
\item [(B)] The axis of $R_0$ is not contained in the mirror of a plane reflection in $G_2$.\\[-.25in]
\end{itemize}
If $G_2$ is dihedral, then Lemma~\ref{r0r1}(b) tells us that $r$ is even and $G_{2}\cong D_{r/2}$. Case (A) can only occur if $G_2$ is dihedral, and then $R_0 R_2$ is the plane reflection with mirror $E_1$ and $\langle R_0, G_2 \rangle$ is isomorphic to $D_{r/2} \times C_2$, with the factor $C_2$ generated by $R_{0}R_{2}$. In dealing with case (A), $R_2$ will always denote the plane reflection in $G_2$ whose mirror contains the axis of $R_0$. On the other hand, in case (B), the subgroup $G_2$ may be cyclic or dihedral. If $G_2$ is dihedral, then the axis of $R_0$ must lie ``halfway'' (in $E_1$) between the plane mirrors of two basic generators of $G_2$ inclined at an angle $2\pi/r$; in fact, by Lemma~\ref{groupstruct}(e), $G_2$ must be invariant under conjugation by $R_0$. In this case, $\langle R_0, G_2 \rangle$ is isomorphic to $D_{r}$ and contains a rotatory reflection of order $r$. 

Recall our notation for the elements of the special group $G_*$ of $G$ (see Section~\ref{thegroup}). In particular, if $R$ is a general element of $G$, then $R'$ is its image in $G_*$. Similarly, if $E$ is a plane, then $E'$ will denote its translate through the origin $o$.

It is convenient to assume that $o$ is the base vertex of $\K$. Then $G_*$ contains $R_1$ as an element and $G_2$ as a subgroup; in fact, $R_1=R_1'$ and $R=R'$ for each $R\in G_2$ (for short,  $G_2=G_2'$), since $R_1$ and $R$ fix the base vertex, $o$. It follows that the vertex-figure group $\langle R_1,G_2\rangle$ of $\K$ at $o$ is a subgroup of $G_*$.

In the present context, $G_*$ is a finite irreducible crystallographic subgroup of ${\mathcal O}(3)$ that contains a rotatory reflection, $R_0' R_1$, whose invariant plane, $E_2'$, is perpendicular to the mirror of a plane reflection, $R_1$, and whose period is just that of $R_{0}R_{1}$. This trivially excludes the groups $[3, 3]^+$ and $[3, 4]^+$ as possibilities for $G_*$ (see (\ref{finthree})). Moreover, we can also eliminate $[3, 3]^*$ on the following grounds. In fact, if we place the vertices of a tetrahedron $\{3, 3\}$ at alternating vertices of a cube that is centered at $o$ and has its faces parallel to the coordinate planes, then the possible mirrors for plane reflections in $[3, 3]^*$ are just the coordinate planes; however, these are not perpendicular to invariant planes of rotatory reflections in $[3, 3]^*$. (Note here that the rotary reflection $-I$, of period $2$, cannot occur 
as $R_{0}'R_{1}$.)  Thus $G_*$ must be one of the groups $[3,3]$ or $[3,4]$.

Next observe that $[3, 3]$ does indeed occur as the special group of an infinite regular polyhedron with finite faces and mirror vector $(1,2)$; however, there is just one polyhedron of this kind, namely $\{6,6\}_4$ (see \cite[p.225]{rap}). On the other hand, the following considerations will exclude $[3, 3]$ as a special group for complexes $\K$ that are not  polyhedra. Suppose $\K$ has $[3,3]$ as its special group but $\K$ is not a polyhedron (that is, $r \geq 3$).

Again it is convenient to place the vertices of $\{3, 3\}$ at alternating vertices of a cube that is centered at $o$ and has its faces parallel to the coordinate planes. Then $R_0'$ is a half-turn about a coordinate axis, $R_0' R_1$ is a rotatory reflection of period $4$ with invariant plane $E_2'$ containing this axis, and $R_1$ is a reflection in a plane perpendicular to $E_2'$. Since $R_0'$ and $R_1$ cannot commute, the mirror of $R_1$ cannot be a coordinate plane.

\begin{figure}
\begin{center}
\includegraphics[width=9cm, height=6cm]{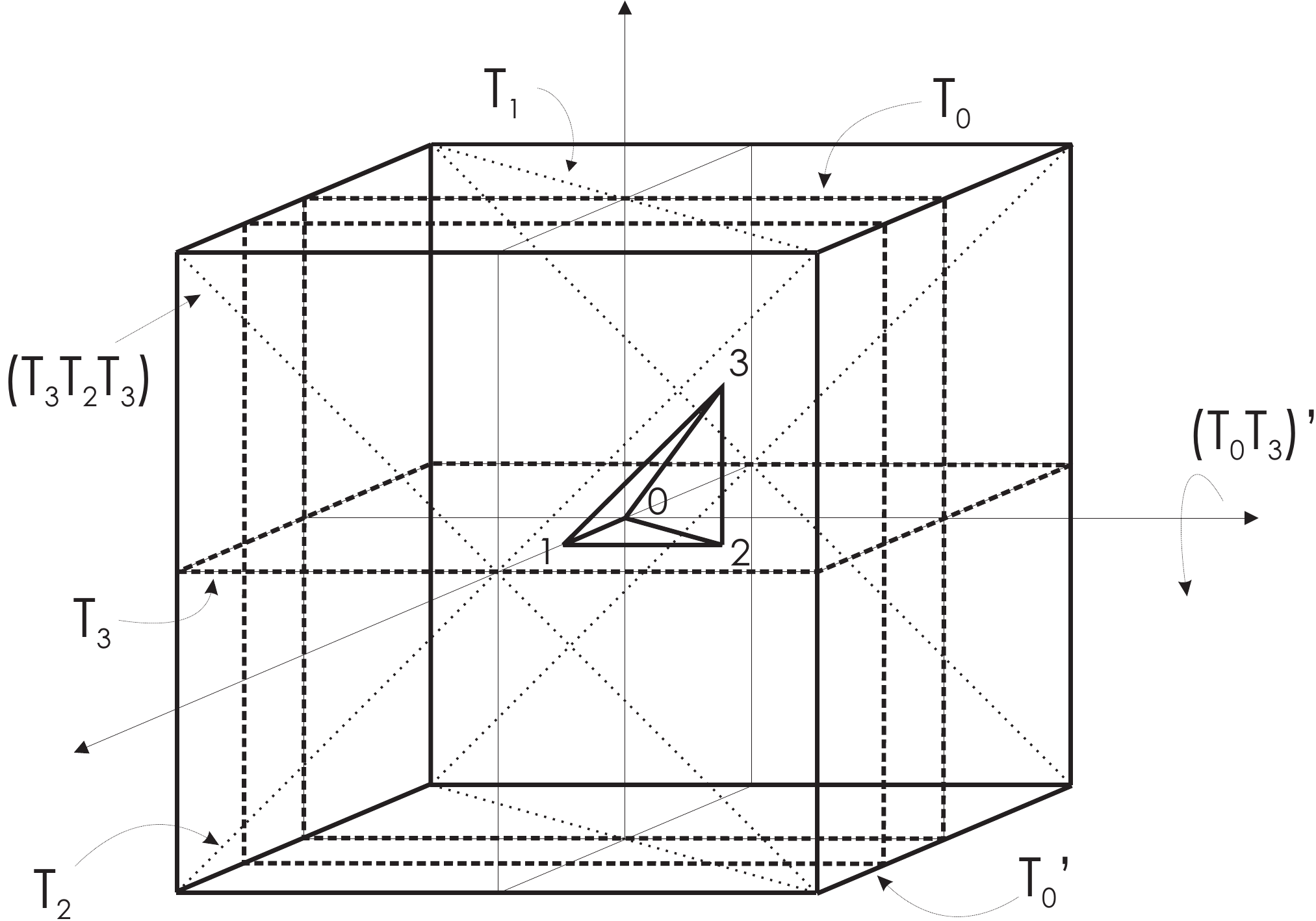}
\caption{Generators of the special group of $H$.}\label{fig8b}
\end{center}
\end{figure}

We now consider the subgroup $G_2$ of $G_*$, bearing in mind that $r\geq 3$. Then $G_2$ must contain a non-trivial rotation $S$ about an axis orthogonal to the axis of $R_0'$. Hence, since $G_*= [3, 3]$, this must necessarily be a half-turn about a coordinate axis distinct from the axis of $R_0'$. Moreover, since $G_*$ is linearly irreducible, we must have $E_1' \ne E_2'$, that is, $E_1'$ is perpendicular to $E_2'$ and the rotation axis of $S$ lies in $E_2'$. Appealing again to irreducibility (or the fact that $\K$ is not a polyhedron), we see that $G_2$ must be a dihedral group of order $4$ and hence contain two plane reflections whose mirrors are perpendicular and intersect in the rotation axis of $S$; moreover, in case (A), the group $G_*$ also contains the reflection in $E_{1}'$. However, since $[3,3]$ does not contain three plane reflections with mutually perpendicular mirrors, this rules out case (A) and only leaves case (B) with a configuration of mirrors and axes as in Figure~\ref{fig8b}, with $R_0'$ and $R_1$ corresponding to $(T_0 T_3)'$ and $T_1$, respectively, and the two reflections in $G_2$ (with mirrors meeting at  the axis of $S$) corresponding to $T_2$ and $T_3 T_2 T_3$. Now, to complete the proof we require the following lemma, which is of interest in its own right.

\begin{lemma}\label{grouph}
Let $[4, 3, 4] = \langle T_0, T_1, T_2, T_3 \rangle$ be the symmetry group of the cubical tessellation $\{4, 3, 4\}$ of $\E$, where $T_0, T_1, T_2, T_3$ are the distinguished generators (see Figure~\ref{fig8b}). Let $H := \langle T_0T_3, T_1, T_2, T_3 T_2 T_3 \rangle$. Then $H$ is a subgroup of $[4, 3, 4]$ of index $2$ and acts simply flag-transitively on the $2$-skeleton of $\{4, 3, 4\}$. In particular, the $2$-skeleton of $\{4,3,4\}$ can be recovered from $H$ by Wythoff's construction.
\end{lemma}

\begin{proof}
Figure~\ref{fig8b} shows the eight cubes of the cubical tessellation $\{4,3,4\}$ that meet at the base vertex $o$ (labeled $0$). The generators $T_0, T_1, T_2, T_3$ and their images $T_0'$, $T_1'=T_1$, $T_2'=T_2$, $T_3'=T_3$ in the special group are indicated. Here, $T_j$ is the reflection in the plane of the small fundamental tetrahedron opposite to the vertex labeled $j$, for $j=0,1,2,3$. 

The two elements $T_1 (T_2 T_3)^2 T_1$ and $T_0 T_3$ are half-turns about parallel axes (in Figure~\ref{fig8b}, the $y$-axis and the line through the points labeled $1, 2$) and their product is a basic translation (along the $x$-axis) in $H$, whose conjugates under $H$ generate the full translation subgroup of $[4,3,4]$. On the other hand, $H$ must be a proper subgroup of $[4,3,4]$; in fact, its special group is a subgroup $[3, 3]$ of the special group $[3, 4]$ of $[4,3,4]$ generated by $(T_0T_3)', T_1, T_2, T_3T_2T_3$. Bearing in mind that $[4,3,4]$ is the semi-direct product of its translation subgroup by its vertex-figure group $[3, 4]$, it now follows that $H$ must have index $2$. The coset of $H$ in $[4,3,4]$ distinct from $H$ is $T_3 H$.

In order to recover the $2$-skeleton of $\{4,3,4\}$ from $H$ define subgroups $H_0 := \langle T_0T_3 \rangle$, $H_1 := \langle T_1 \rangle$ and $H_2 := \langle T_2, T_3 T_2 T_3 \rangle$. Then the base vertex $o$ of $\{4,3,4\}$ is fixed by $H_1$ and $H_2$, and the base edge of $\{4,3,4\}$ (through the points labeled $0$ and $1$) is fixed by $H_2$. Applying Wythoff's construction with group $H$ (with distinguished subgroups $H_{0}$, $H_{1}$, $H_2$ as indicated) and with initial vertex $o$, then yields a complex with its vertices, edges and $2$-faces among those of $\{4, 3, 4\}$. However, since $T_3 H$ is the only non-trivial coset of $H$ and the element $T_3$ leaves the base vertex, base edge and base $2$-face of $\{4,3,4\}$ invariant, the images of the latter under the full group $[4,3,4]$ are just those under the subgroup $H$. Hence, Wythoff's construction applied to $H$ gives the full $2$-skeleton of $\{4,3,4\}$. Moreover, since $T_3$ is the reflection in the plane through the base $2$-face of $\{4,3,4\}$, it stabilizes the base flag of the $2$-skeleton and hence is the only non-trivial element in the flag stabilizer. Since $T_3$ is not in $H$, it follows that $H$ must be simply flag-transitive.
\end{proof} 

We now complete the proof that $[3,3]$ cannot occur as a special group $G_*=\langle R_{0}',R_{1},G_{2}\rangle$ if $\K$ is not a polyhedron. Recall that we were in case (B). As noted earlier, the generators $R_0', R_1$ and the two reflections in $G_2$ can be viewed as appropriate elements of the special group $H_*$ of the simply flag-transitive subgroup $H$ of the symmetry group of the $2$-skeleton of $\{4,3,4\}$ described in Lemma \ref{grouph} and depicted in Figure~\ref{fig8b}. In particular, we may identify the base vertex and base edge of $\K$ with those of the $2$-skeleton of $\{4,3,4\}$; recall here that the axis of $S$ is simply the line through the base edge of $\{4,3,4\}$. A straightforward application of Wythoff's construction then shows that $\K$ must indeed coincide with the $2$-skeleton of $\{4,3,4\}$, contradicting our basic assumption that $\K$ is simply flag-transitive (recall that the latter means that the {\em full\/} symmetry group is simply flag-transitive).

In summary, we have established the following

\begin{lemma}\label{finfacspec}
Let $\K$ be a simply flag-transitive regular complex with finite faces and mirror vector $(1,2)$, and let $\K$ not be the polyhedron $\{6,6\}_4$. Then $G_{*}=[3,4]$. 
\end{lemma}

\subsection{The complexes associated with $[3, 4]$}
\label{com34}

As before, let $\K$ be a simply flag-transitive regular complex $\K$ with finite faces and mirror vector $(1,2)$. Recall our standing assumption that $\K$ is infinite and $G$ is irreducible. Suppose the special group $G_*$ of $\K$ is the octahedral group $[3,4]$. From the previous lemma we know that this only eliminates the polyhedron $\{6,6\}_4$. We further assume that the base vertex $F_0$ of $\K$ is the origin $o$, so in particular $R_1'=R_1$ and $G_2'=G_2$. The vertex $v$ of the base edge $F_1$ of $\K$ distinct from $F_0$ is called the {\em twin vertex\/} of $\K$ (with respect to the base flag).

First observe that $G_2$ must contain a non-trivial rotation. Otherwise, $G_2$ is generated by a single plane reflection and so has order $2$. Furthermore, $\K$ must then be a (pure) regular polyhedron with complete dimension vector $(1,2,2)$ (see \cite[p.225]{rap}); however, no such polyhedron exists. Thus $G_2$ contains a rotation, $S$ (say), that generates its non-trivial rotation subgroup. Note that, as a vector, the twin vertex spans the rotation axis of $S$ passing through $o$ and $v$.

As the reference figure for the action of the special group $G_*$ we take the cube $C:=\{4,3\}$ with vertices $(\pm 1,\pm 1,\pm 1)$, so that $C$ is centered at $o$ and has faces parallel to the coordinate planes. Its group $[3, 4]$ contains six rotatory reflections of period $4$ and eight rotatory reflections of period $6$. Each rotatory reflection of period $4$ is given by a rotation by $\pm \pi/2$ about a coordinate axis, followed by a reflection in the coordinate plane perpendicular to the axis. Each rotatory reflection of period $6$ is given by a rotation by $\pm \pi/3$ about a main diagonal of $C$, followed by a reflection in the plane through $o$ perpendicular to the diagonal.

We break the discussion down into two cases, Case I and Case II, according as $R_0' R_1$ is a rotatory reflection of period $4$ or $6$, or, equivalently, $\K$ has square faces or hexagonal faces. (It is convenient here to use the term ``square" to describe a regular quadrangle, such as the Petrie polygon of $\{3,3\}$.)  Recall that $E_1'$ denotes the plane through $o$ parallel to the perpendicular bisector $E_1$ of $F_1$ (on which $G_2$ acts faithfully), and that $E_2'$ is the invariant plane of $R_0' R_1$ and contains the axis $E_1' \cap E_2'$ of the half-turn $R_0'$.

\medskip
\noindent{\bf Case I: Square faces}
\medskip

Suppose $R_0' R_1$ is a rotatory reflection of period $4$ with invariant plane $E_2'$ given by the $xy$-plane (say). Then, up to conjugacy in $[3,4]$, there are exactly two possible choices for the half-turn $R_0'$, namely $R_0'$ either rotates about the center of a face of $C$ (a coordinate axis) or about the midpoint of an edge of $C$. 

For the sake of simplicity, when we claim uniqueness for the choice of certain elements within the special group or of mirrors or fixed point sets of such elements, we will usually omit any qualifying statements such as ``up to conjugacy" or ``up to congruence". Throughout, these qualifications are understood. 
 
\medskip
\noindent{\em Case Ia:  $R_0'$ rotates about the center of a face of $C$}
\medskip

Suppose $R_{0}'$ rotates about the $y$-axis (say). Since the mirror of $R_1$ is a plane perpendicular to $E_2'$ and the elements $R_0'$ and $R_1$ do not commute, there is only one choice for $R_1$, namely the reflection in the plane $y=x$. We next consider the possible choices for $S$ and $E_1'$. Clearly, $E_1'\neq E_2'$, since otherwise $E_2'$ is an invariant subspace for the (irreducible) group $G_{*}$. Bearing in mind that $E_1'$ must contain the rotation axis of $R_0'$, this leaves the following two possibilities for $E_1'$.

First, suppose $E_1'$ is perpendicular to $E_2'$. Then $E_1'$ must be the $yz$-plane. We can eliminate this possibility as follows. 

Suppose $G_2$ is cyclic of order $r$. Then $r\neq 2$, once more since otherwise $E_2'$ is an invariant subspace of $G_*$. If $r>2$, then necessarily $r=4$ and the twin vertex $v$ (invariant under $S$) is of the form $v=(a,0,0)$ for some non-zero parameter $a$. In particular, this forces the square faces of $\K$ to be planar and parallel to the coordinate planes. In fact, Wythoff's construction shows in this case that $\K$ must be the $2$-skeleton of the cubical tessellation $\{4,3,4\}$ with vertex-set $a\mathbb{Z}^3$. Moreover, $G$ must be the full group $[4,3,4]$, since its full translation subgroup and special group are the same as those of $\{4,3,4\}$; note here that the translation by $v$ belongs to $G$ (the axes of the half-turns $R_1S^2R_1$ and $R_0$, respectively, are the lines through $o$ and $\frac{1}{2}v$ parallel to the $y$-axis, so $R_1S^2R_1R_0$ is the translation by $v$).  All this contradicts our basic assumption that $\K$ is simply flag-transitive.

Next consider the possibility that $G_2$ is dihedral of order $r$. Then $r=4$ or $8$ and we are in case (A) or (B), depending on whether or not the axis of $R_0$ lies in the mirror of a plane reflection of $G_2$ (see Section~\ref{specgroup}). In case (A) we cannot have $r=4$, once again on account of the irreducibility of $G_*$. In case (A) with $r=8$, the resulting group $G$ is a supergroup of the group discussed in the previous paragraph. In particular, the twin vertex is of the form $v=(a,0,0)$ and $\K$ has face mirrors; for example, $F_2$ lies in the $xy$-plane, which is the mirror of the reflection $R_2$ in $G_2$. Thus $\K$ is the $2$-skeleton of $\{4,3,4\}$ and $G=[4,3,4]$, again a contradiction to our assumptions on $\K$. In any case, case (A) with $r=8$ can be ruled out and we are left with case (B) with  $r=4$. There is just one possible configuration, namely a pair of generating reflections of $G_2$ with mirrors given by the planes $z=y$ and $z=-y$. This is exactly the mirror configuration described in Lemma~\ref{grouph} and depicted in Figure~\ref{fig8b}. However, now $R_0'$, $R_1$ and $G_2$ all preserve the two sets of alternating vertices of $C$, so $G_*=[3,3]$, not $[3,4]$, allowing us to eliminate case (B) as well.

Second, suppose $E_1'$ is not perpendicular to $E_2'$. Then $E_1'$ is the plane $z=-x$ (say) and $S$ is a half-turn about the midpoint $(1,0,1)$ of an edge of $C$. In any case, $r=2$ or $4$ and the twin vertex has the form $v=(a,0,a)$ with $a\neq 0$. 

If $r=2$ (and hence $G_2$ is cyclic), the complex $\K$ is a polyhedron with skew square faces and planar hexagonal vertex-figures. In particular, $\K$ is the polyhedron $\{4,6\}_6$, since its complete dimension vector is $(1,2,1)$ and its special group is $[3,4]$ (see \cite[p.225]{rap}). Here the square faces of $\K$ are inscribed into three quarters of all cubes of a cubical tessellation in $\E$. 

\begin{figure}
\begin{center}
\includegraphics[width=9cm, height=7cm]{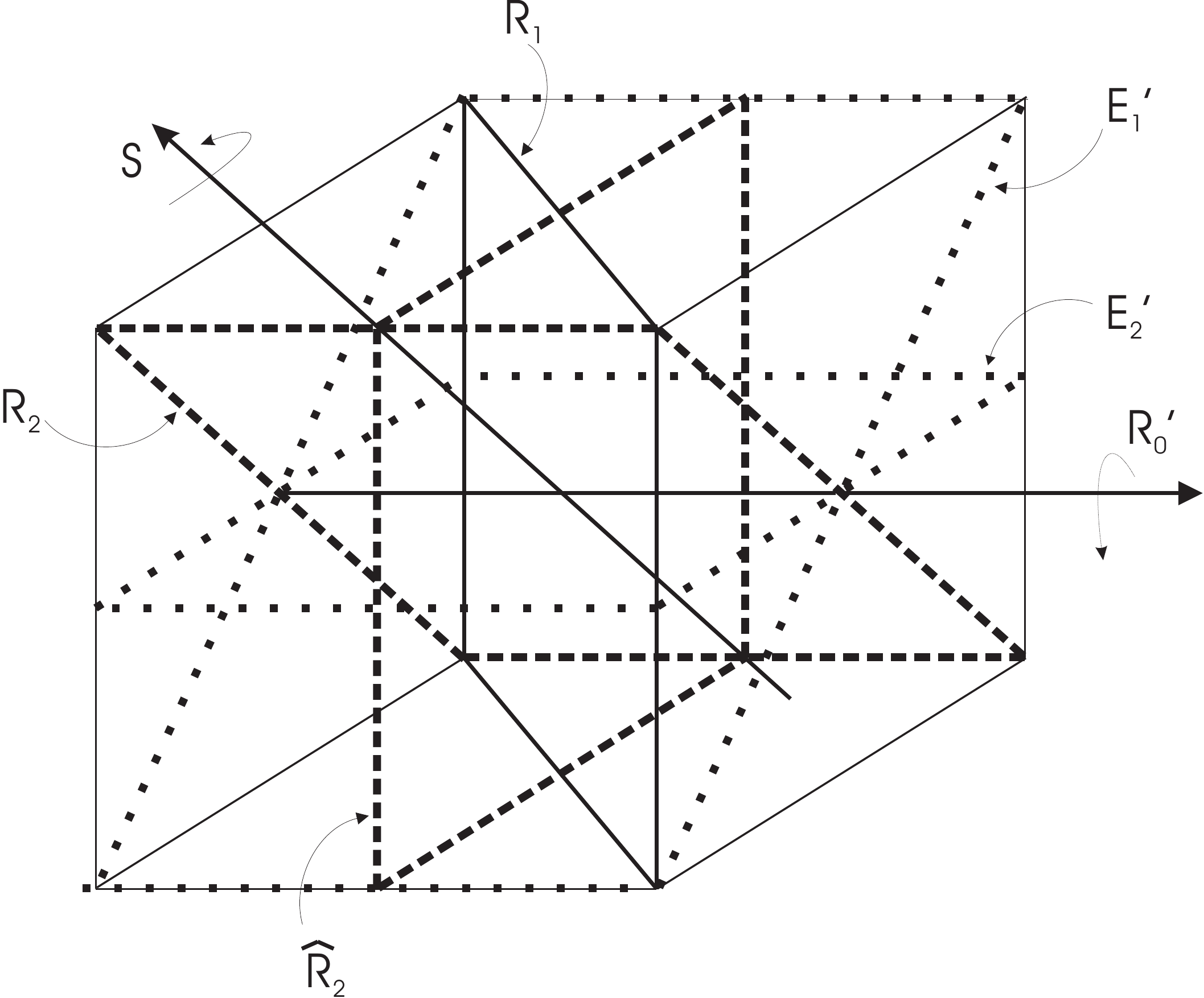}
\caption{The special group of the complex $\K_{1}(1, 2)$}\label{k112}
\end{center}
\end{figure}

If $r=4$, then $G_2$ is dihedral and case (B) cannot occur.  On the other hand, in case (A), the group $G_2$ is generated by the reflections $R_2,\hat{R}_2$ in the planes $z=x$ and $y=0$, respectively, and the mirror of $R_2$ contains the axis of $R_0'$ (see Figure \ref{k112}). Then the set of generators $R_0,R_1,R_2,\hat{R}_2$ of $G$ is given by 
\begin{equation}
\label{1gen12}
\begin{array}{rccl}
R_0\colon    & (x,y,z)  &\mapsto & (-x,y,-z) + (a,0,a),\\
R_{1}\colon & (x,y,z)  &\mapsto  & (y,x,z),\\
R_{2}\colon & (x,y,z)   &\mapsto & (z,y,x),\\
\hat{R}_2\colon & (x,y,z) & \mapsto & (x,-y,z),  \\
\end{array}
\end{equation}
with $a\neq 0$. Now we obtain a legitimate regular complex, denoted $\K_{1}(1,2)$, with skew square faces, four around each edge. (Recall our convention to label regular complexes with their mirror vectors.)

The vertex-set of $\K_{1}(1,2)$ consists of the images of the base vertex $o$ under $G$ and coincides with the lattice $\Lambda_{(a,a,0)}$. The faces of $\K_{1}(1,2)$ are Petrie polygons of tetrahedra inscribed at alternating vertices of  cubes in the cubical tessellation with vertex-set $a\mathbb{Z}^3$. In each cube, $\K_{1}(1,2)$ takes edges of just one of the two possible tetrahedra; the Petrie polygons of this tetrahedron determine three faces of $\K_{1}(1,2)$, forming a finite subcomplex $\{4,3\}_3$ (the Petrie dual of $\{3,3\}$), such that every edge lies in exactly two of them. Thus every edge of $\K_{1}(1,2)$ is surrounded by four faces, so that each cube in a pair of adjacent cubes contributes exactly two faces to the unique edge of $\K_{1}(1,2)$ contained in their intersection. In particular, the vertices of the base face $F_2$ are the images of $o$ under $\langle R_0,R_1\rangle$ and are given by $o$, $v=(a,0,a)$, $(a,a,0)$, $(0,a,a)$, occurring in this order. Moreover, the vertex-figure group $\langle R_1,G_2\rangle$ of $\K_{1}(1,2)$ at $o$ is the octahedral group $[3,4]$, occurring here with (standard) generators $R_2,R_1,\hat{R}_2$. The vertex-figure of $\K_{1}(1,2)$ at $o$ is the (simple) edge graph of the cuboctahedron with vertices $(\pm a,\pm a, 0)$, $(0,\pm a,\pm a)$, $(\pm a,0,\pm a)$; in fact, each edge of this graph corresponds to a face of $\K_{1}(1,2)$ with vertex $o$, and vice versa.

\medskip
\noindent{\em Case Ib:  $R_0'$ rotates about the midpoint of an edge of $C$}
\medskip

Recall that the rotatory reflection $R_0' R_1$ of period $4$ has the $xy$-plane as its invariant plane $E_2'$. Here we may assume that $R_0'$ is the half-turn about the line through $o$ and the midpoint $(1,1,0)$ of an edge of $C$, and that $R_1$ is the reflection in a coordinate plane perpendicular to $E_2'$, the plane $y=0$ (say). Since $E_1'$ must contain the axis of $R_0'$, there are two possible choices for $E_1'$, namely either $E_1'$ is perpendicular to $E_2'$ or to a main diagonal of $C$; note here that we must have $E_1'\neq E_2'$, since $G$ is irreducible. We show that the first possibility cannot occur and that the second contributes a new regular complex.

First, suppose $E_1'$ is perpendicular to $E_2'$. Then $E_1'$ is the plane $y=x$ (say), and $G_2$ is either cyclic of order $2$ or dihedral of order $4$. In either case, $G_*$ would be reducible, with invariant subspace $E_2'$. Thus $E_1'$ cannot be perpendicular to $E_2'$.

Second, suppose $E_2'$ is perpendicular to a main diagonal of $C$, the diagonal passing through the pair of antipodal vertices $\pm (1,-1,1)$ (say). Then $r=3$ or $6$, and $G_2$ is either cyclic of order $3$ or dihedral of order $6$, the latter necessarily occurring here as case (B). In any case, the twin vertex of $\K$ is of the form $v=(a,-a,a)$ with $a\neq 0$. 

\begin{figure}
\begin{center}
\includegraphics[width=8cm, height=6cm]{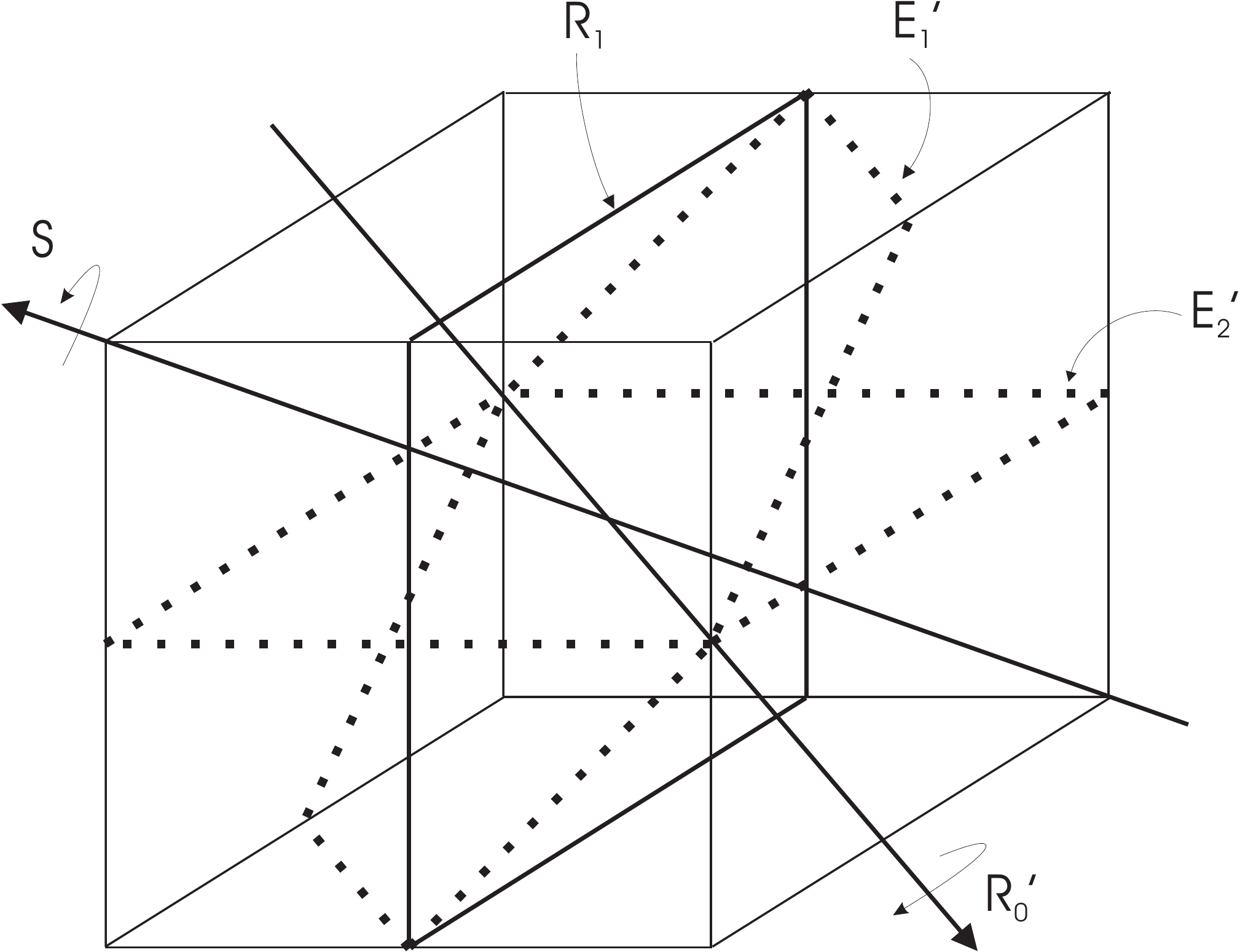}
\caption{The special group of the complex $\K_{2}(1, 2)$}\label{k212}
\end{center}
\end{figure}

Suppose $G_2$ is cyclic of order $3$ and is generated by a $3$-fold rotation $S$ about the main diagonal of $C$ passing through $\pm (1,-1,1)$. Then the generators $R_0,R_1,S$ of $G$ are given by
\begin{equation}
\label{2gen12}
\begin{array}{rccl}
R_0\colon    & (x,y,z)  &\mapsto & (y,x,-z) + (a,-a,a),\\
R_{1}\colon & (x,y,z)  &\mapsto  & (x,-y,z),\\
S\colon  & (x,y,z) &\mapsto & (-y,-z,x),
\end{array}
\end{equation}
with $a\neq 0$ (see Figure \ref{k212}). This determines a new regular complex, denoted $\K_{2}(1,2)$, with skew square faces, three around each edge.

The vertex-set of $\K_{2}(1,2)$ is the lattice $\Lambda_{(a,a,a)}$. Relative to the cubical tessellation of $\E$ with vertex set $(a,a,a)+2a\mathbb{Z}^3$, a typical square face of $\K_{2}(1,2)$ has its four vertices, in pairs of opposites, located at the centers and at common neighboring vertices of adjacent cubes. In this way, each pair of adjacent cubes determines exactly two faces of $\K_{2}(1,2)$, and these meet only at their two vertices located at the centers of the cubes. Moreover, the three cubes that share a common vertex with a given cube and are adjacent to it, determine the three faces of $\K_{2}(1,2)$ that surround the edge that joins the common vertex to the center of the given cube. In particular, the vertices of $F_2$ are $o$, $v=(a,-a,a)$, $(2a,0,0)$, $(a,a,a)$, in this order. Now the vertex-figure group $\langle R_1,G_2\rangle$ of $\K_{2}(1,2)$ at $o$ is the group $[3, 3]^*$ of order $24$. The vertex-figure at $o$ itself is the (simple) edge graph of the cube with vertices $(\pm a,\pm a, \pm a)$.

It is interesting to observe the following nice picture of the face of $\K_{2}(1,2)$. The fundamental tetrahedron of the Coxeter group $P_4$ whose diagram is an unmarked square has two opposite edges of length $2$ and four of length $\sqrt{3}$ (when defined relative to $2\mathbb{Z}^3$). Up to scaling, these four others give the shape of the face of $\K_{2}(1,2)$. The same remark applies to the next complex, $\K_3(1,2)$.

\begin{figure}
\begin{center}
\includegraphics[width=9cm, height=7cm]{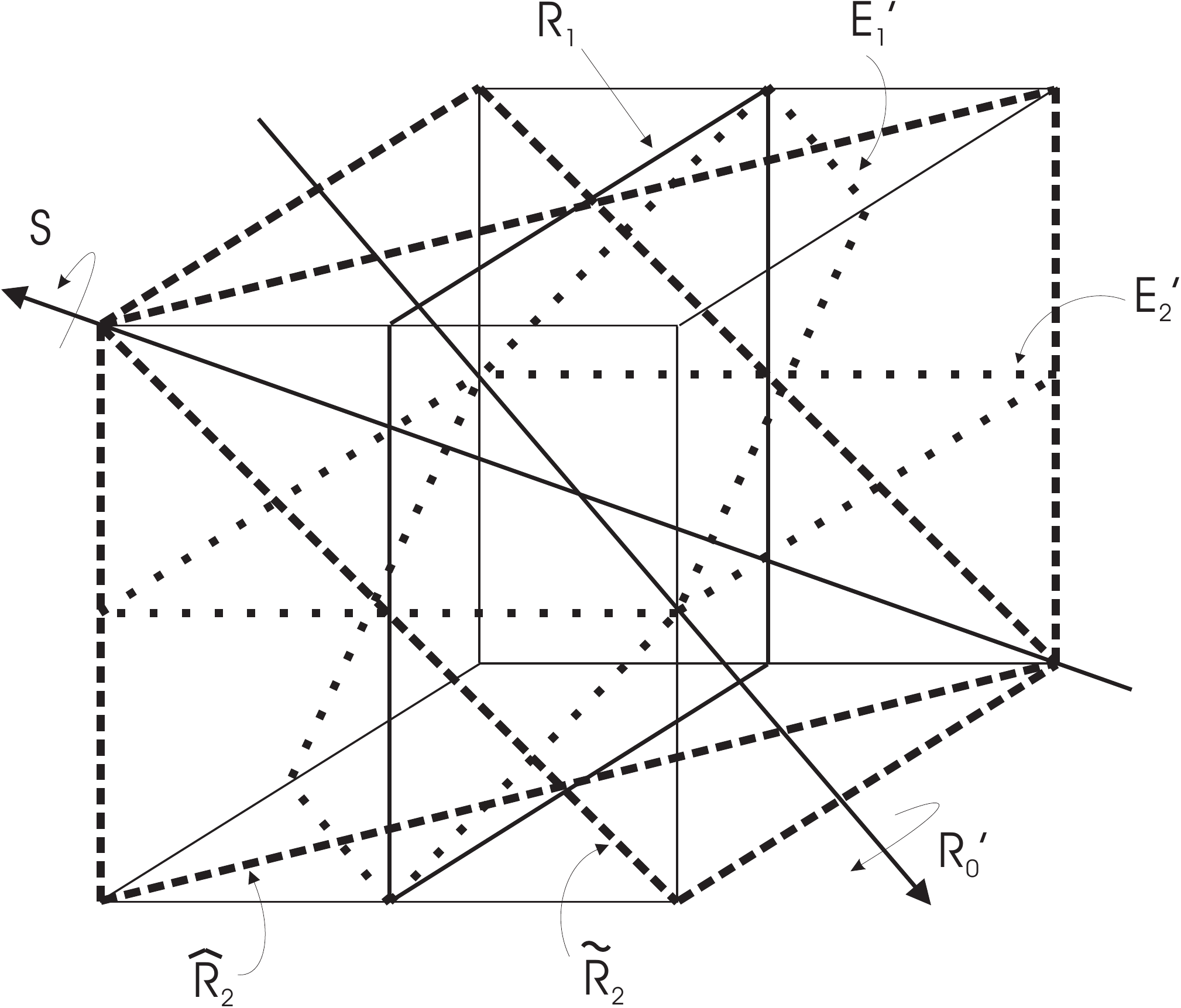}
\caption{The special group of the complex $\K_{3}(1, 2)$}\label{k312}
\end{center}
\end{figure}

Finally, let $G_2$ be dihedral. As generators of $G_2$ we take the reflections $\hat{R}_{2}$ and $\widetilde{R}_2$ in the planes $y=-x$ and $z=-y$, respectively, so that $\hat{R}_{2}\widetilde{R}_2=S$ with $S$ as in (\ref{2gen12}). Then $G$ has generators $R_0,R_1,\hat{R}_{2},\widetilde{R}_2$ given by
\begin{equation}
\label{3gen12}
\begin{array}{rccl}
R_0\colon    & (x,y,z)  &\mapsto & (y,x,-z) + (a,-a,a),\\
R_{1}\colon & (x,y,z)  &\mapsto  & (x,-y,z),\\
\hat{R}_{2}\colon  & (x,y,z) &\mapsto & (-y,-x,z),\\
\widetilde{R}_2\colon & (x,y,z)  &\mapsto & (x,-z,-y),
\end{array}
\end{equation}
again with $a\neq 0$ (see Figure \ref{k312}). They also determine a regular complex, denoted $\K_{3}(1,2)$, with skew square faces, now six around each edge. This contains $\K_{2}(1,2)$ as a subcomplex.

As for the previous complex, the vertex-set of $\K_{3}(1,2)$ is $\Lambda_{(a,a,a)}$ and a typical square face has its four vertices, in pairs of opposites, located at the centers and at common neighboring vertices of adjacent cubes in the cubical tessellation with vertex set $(a,a,a)+2a\mathbb{Z}^3$. Now each pair of adjacent cubes determines four faces of $\K_{3}(1,2)$, not only two as for $\K_{2}(1,2)$. The three cubes that share a vertex with a given cube and are adjacent to it, now determine all six faces of $\K_{3}(1,2)$ that surround the edge that joins the common vertex to the center of the given cube. The base face $F_2$ is the same as for $\K_{2}(1,2)$. However, the vertex-figure group $\langle R_1,G_2\rangle$ of $\K_{3}(1,2)$ at $o$ now is the full group $[3,4]$, so the vertex-figure at $o$ itself is the {\em double-edge graph\/} of the cube with vertices $(\pm a,\pm a, \pm a)$, meaning the graph obtained from the ordinary edge graph by doubling the edges but maintaining the vertices. (Note that the double-edge graph admits an action of the vertex-figure group.)

Observe also that the faces of $\mathcal K_3(1,2)$ are those of the Petrie duals of the facets $\{\infty, 4\}_4 \# \{\, \}$ of the $4$-apeirotope $\apeir\{4,3\}$, with $\{4,3\}$ properly chosen.

\medskip
\noindent{\bf Case II: Hexagonal faces}
\medskip

Suppose $R_0' R_1$ is a rotatory reflection of period $6$ whose invariant plane $E_2'$ is perpendicular to a main diagonal of $C$, the diagonal passing through the vertices $\pm (1,1,1)$ (say). Then we may take $R_0'$ to be the half-turn about the line through $o$ and $(1,0,-1)$ (say), the latter being the midpoint of an edge of $C$. Moreover, since the mirror of $R_1$ is perpendicular to $E_2'$ and $R_0' R_1$ has period $6$, we may assume $R_1$ to be the reflection in the plane $y=x$. With $R_0',R_1$ specified, we now have three choices for $E_1'$, namely $E_1'$ is a coordinate plane or is perpendicular to either a main diagonal of $C$ or to a line through the midpoints of a pair of antipodal edges of $C$.

\medskip
\noindent{\em Case IIa: $E_1'$ is a coordinate plane}
\medskip

Suppose $E_1'$ is the $xz$-plane (say), so that $S$ rotates about the $y$-axis. Then $G_2$ is cyclic of order $2$ or $4$, or dihedral of order $4$ or $8$. In either case, the twin vertex is of the form $v=(0,a,0)$ with $a\neq 0$.

If $G_2$ is cyclic of order $2$, then $\K$ is a polyhedron with complete dimension vector $(1,2,1)$. In fact, comparison with \cite[p.225]{rap} shows that $K = \{6,4\}_6$. The faces of $\K$ are skew hexagons given by Petrie polygons of half the cubes in a cubical tessellation, and the vertex-figures are planar squares.

Next we eliminate the possibility that $G_2$ is cyclic of order $4$. In fact, the presence of $S$ as a symmetry of $\K$ already forces $G_2$ to be dihedral; that is, $G_2$ cannot be cyclic.  In the present configuration, $S$ and $R_1$ already generate the maximum possible group, namely the full octahedral group $[3,4]=G_*$ (see Lemma \ref{casesnum1}). Since the vertex-figure group  $\langle R_1,G_2\rangle$ of $\K$ at $o$ must be a subgroup of $G_*$ containing $S$, it must coincide with $G_*$ and hence contain the full dihedral subgroup $D_4$ of $G_*$ that contains $S$. Thus $D_4$ is a subgroup of $G$. Moreover, since $D_4$ consists of symmetries of $\K$ that fix $o$ and $v$, it lies in the pointwise stabilizer $G_2$ of $F_1$. Hence $G_2=D_4$.

\begin{figure}
\begin{center}
\includegraphics[width=8.5cm, height=7cm]{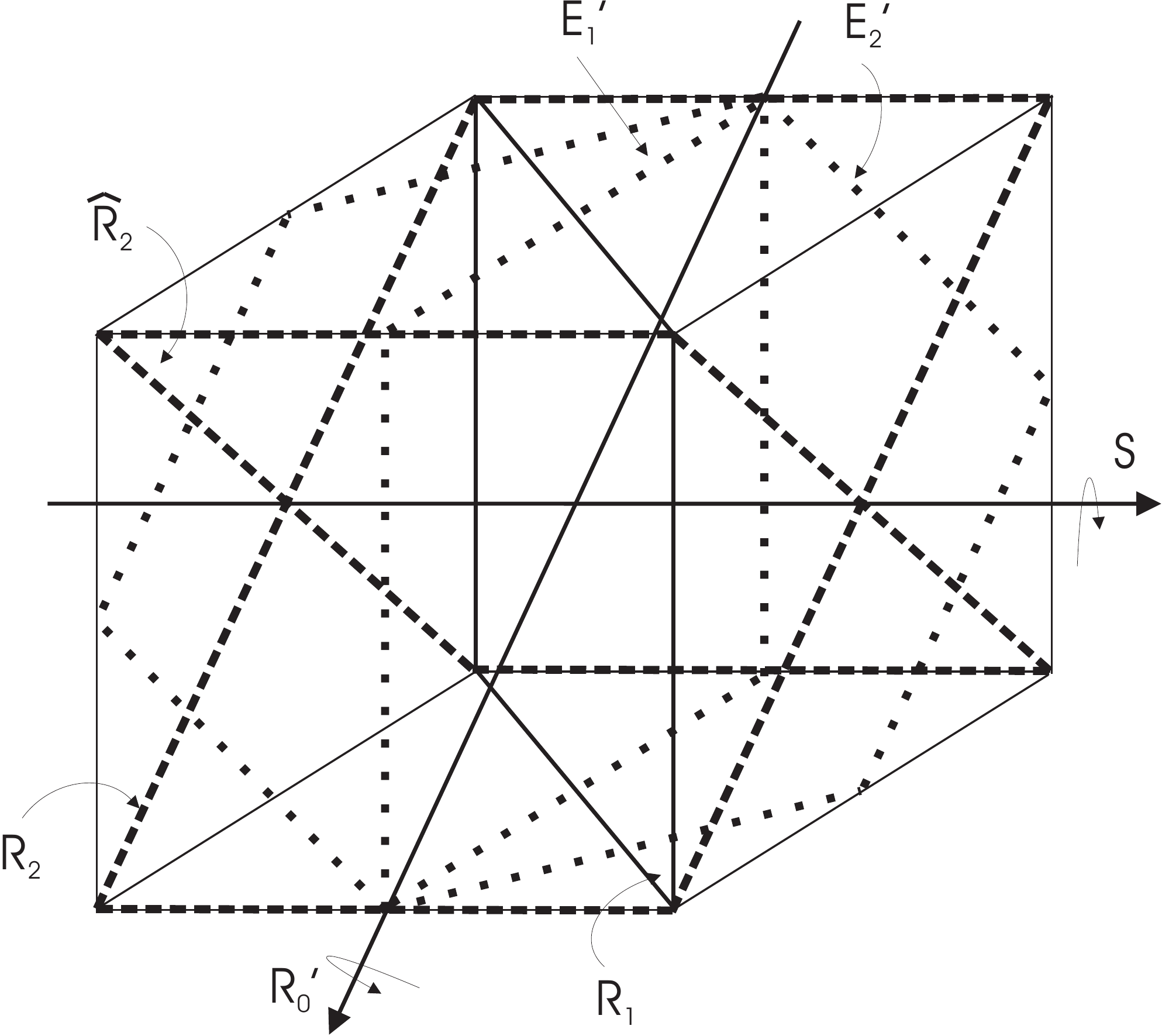}
\caption{The special group of the complex $\K_{4}(1, 2)$}\label{k412}
\end{center}
\end{figure}

Now suppose that $G_2$ is dihedral of order $4$. Then, in case (A), $G_2$ is generated by the reflections $R_2$ and $\hat{R}_{2}$ in the planes $z=-x$ and $z=x$, respectively. Hence the generators $R_0,R_1,R_2,\hat{R}_{2}$ of $G$ are given by
\begin{equation}
\label{4gen12}
\begin{array}{rccl}
R_0\colon    & (x,y,z)  &\mapsto & (-z,-y,-x) + (0,a,0),\\
R_{1}\colon & (x,y,z)  &\mapsto  & (y,x,z),\\
R_{2}\colon  & (x,y,z) &\mapsto & (-z,y,-x),\\
\hat{R}_2\colon & (x,y,z)  &\mapsto & (z,y,x),
\end{array}
\end{equation}
with $a\neq 0$ (see Figure \ref{k412}). In particular, they yield a regular complex, denoted $\K_{4}(1,2)$, which has skew hexagonal faces, four surrounding each edge.

The vertex-set of $\K_{4}(1,2)$ is $a\mathbb{Z}^3$. The faces of $\K_{4}(1,2)$ are Petrie polygons of half the cubes in the cubical tessellation of $\E$ with vertex set $a\mathbb{Z}^3$. In each cube occupied, the Petrie polygons of this cube determine four faces of $\K_{4}(1,2)$, forming a finite subcomplex $\{6,3\}_4$ (the Petrie dual of $\{4,3\}$), such that every edge lies in exactly two of them. Thus the edge graph of $\K_{4}(1,2)$ is the full edge graph of the cubical tessellation and every edge is surrounded by four faces, such that occupied cubes with a common edge (they are necessarily non-adjacent) contribute exactly two faces to this edge of $\K_{4}(1,2)$. In particular, the vertices of the base face $F_2$ are $o$, $v=(0,a,0)$, $(0,a,-a)$, $(a,a,-a)$, $(a,0,-a)$, $(a,0,0)$, in this order. Now the vertex-figure group $\langle R_1,G_2\rangle$ of $\K_{4}(1,2)$ at $o$ is $[3,3]$, occurring here with standard generators $R_2,R_1,\hat{R}_2$. The vertex-figure at $o$ itself is the (simple) edge graph of the octahedron with vertices $(\pm a,0,0)$, $(0,\pm a,0)$, $(0,0,\pm a)$. 

\begin{figure}
\begin{center}
\includegraphics[width=8.5cm, height=7cm]{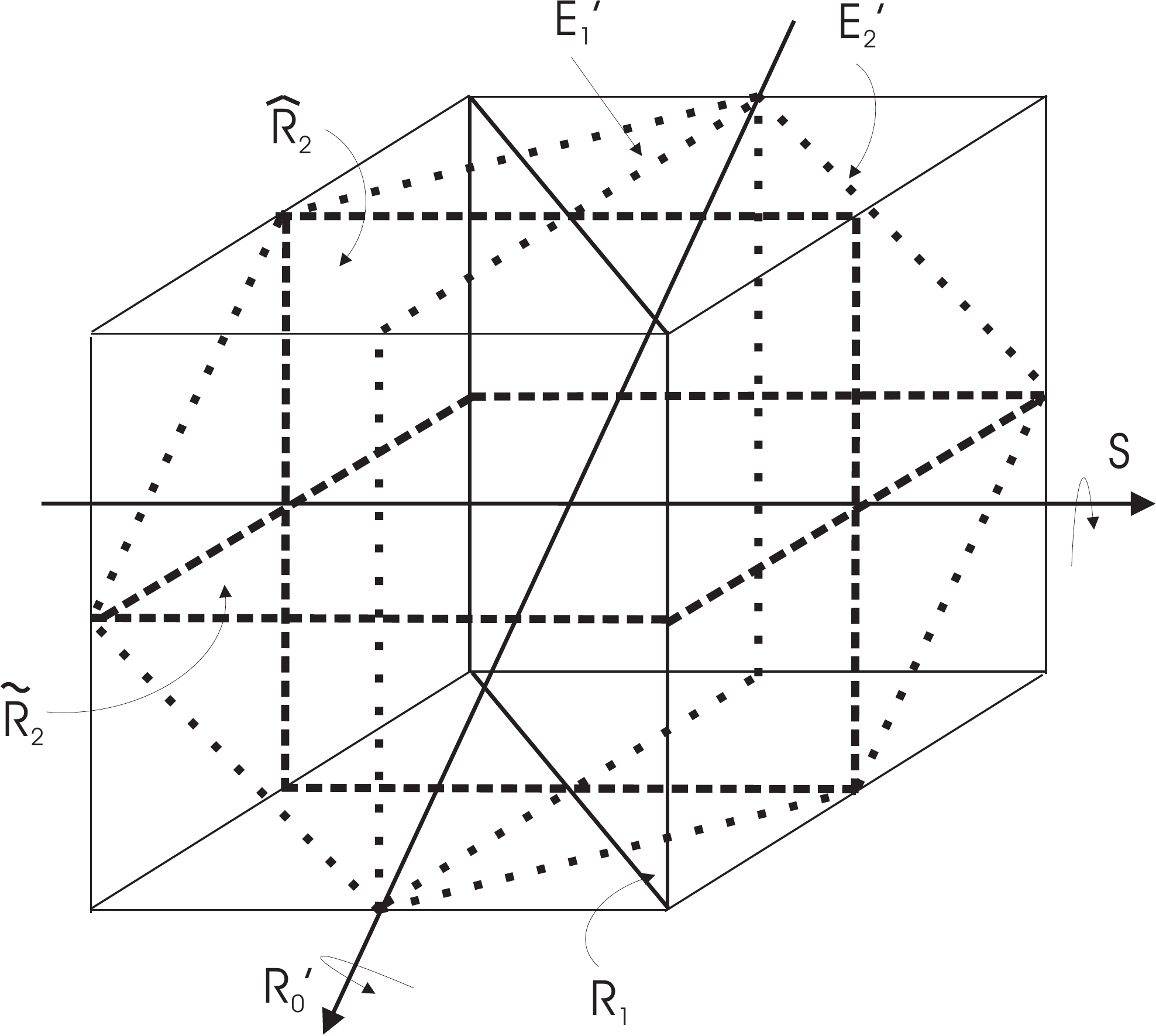}
\caption{The special group of the complex $\K_{5}(1, 2)$}\label{k512}
\end{center}
\end{figure}

Let $G_2$ be dihedral of order $4$ and consider case (B). As generators of $G_2$ we have the reflections ${\hat R}_2,{\widetilde R}_{2}$ in the planes $x=0$ and $z=0$, respectively. Then $R_0,R_1,{\hat R}_2,{\widetilde R}_{2}$ generate $G$ and are given by
\begin{equation}
\label{5gen12}
\begin{array}{rccl}
R_0\colon    & (x,y,z)  &\mapsto & (-z,-y,-x) + (0,a,0),\\
R_{1}\colon & (x,y,z)  &\mapsto  & (y,x,z),\\
{\hat R}_{2}\colon  & (x,y,z) &\mapsto & (-x,y,z),\\
{\widetilde R}_2\colon & (x,y,z)  &\mapsto & (x,y,-z),
\end{array}
\end{equation}
again with $a\neq 0$ (see Figure \ref{k512}). These generators give a regular complex, denoted $\K_{5}(1,2)$, which also has skew hexagonal faces such that four surround each edge.

Now the vertex set of the complex $\K_{5}(1,2)$ is 
$a\mathbb{Z}^3\! \setminus\! \big( (0,0,a) + a\Lambda_{(1,1,1)}\big)$.
As in the previous case, the faces are Petrie polygons of cubes in the cubical tessellation of $\E$ with vertex set $a\mathbb{Z}^3$, with $F_2$ exactly as before. Now each cube occurs and contributes a single face (its Petrie polygon with vertices not in $(0, 0, a) + a\Lambda_{(1, 1, 1)}$), and each vertex lies in exactly eight faces. In particular, the eight faces with vertex $o$ are those Petrie polygons of the eight cubes with vertex $o$ that have $o$ as a vertex and have their two edges incident with $o$ lying in the plane $z=0$. The vertex-figure at $o$ is the double-edge graph of the square with vertices $(\pm a,0,0)$ and $(0,\pm a,0)$, meaning again the ordinary edge graph with its edges doubled. Thus the vertex-figure at $o$ is planar, lying in the plane $z=0$. Observe that the vertex-figure group at $o$ is a reducible group $[4,2] \cong D_{4}\times C_{2}$, occurring here with generators ${\widetilde R}_{2},R_1,{\hat R}_2$ and leaving the plane $z=0$ invariant.

\begin{figure}
\begin{center}
\includegraphics[width=8.5cm, height=7cm]{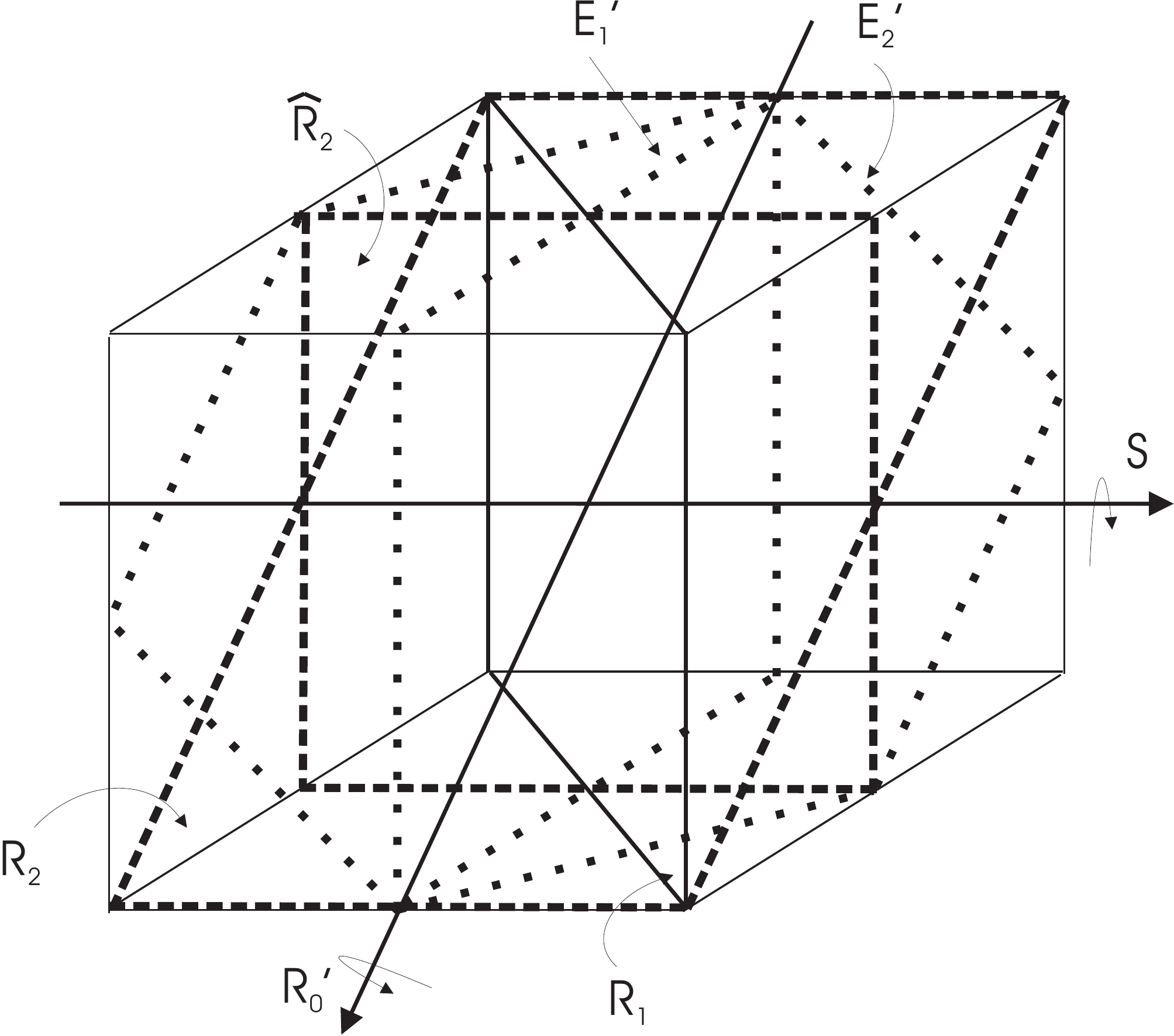}
\caption{The special group of the complex $\K_{6}(1, 2)$}\label{k612}
\end{center}
\end{figure}

It remains to consider the case that $G_2$ is dihedral of order $8$. Then we are in case (A) and $G_2$ is generated by the reflections $R_2$ and $\hat{R}_{2}$ in the planes $z=-x$ and $x=0$, respectively. Now the generators $R_0,R_1,R_2,\hat{R}_{2}$ of $G$ are given by
\begin{equation}
\label{6gen12}
\begin{array}{rccl}
R_0\colon    & (x,y,z)  &\mapsto & (-z,-y,-x) + (0,a,0),\\
R_{1}\colon & (x,y,z)  &\mapsto  & (y,x,z),\\
R_{2}\colon  & (x,y,z) &\mapsto & (-z,y,-x),\\
\hat{R}_2\colon & (x,y,z)  &\mapsto & (-x,y,z),
\end{array}
\end{equation}
again with $a\neq 0$ (see Figure \ref{k612}). Once more we obtain a regular complex, denoted $\K_{6}(1,2)$, again with skew hexagonal faces but now eight surrounding each edge.

The vertex-set of $\K_{6}(1,2)$ is the full lattice $a\mathbb{Z}^3$. As for the two previous complexes, the faces are Petrie polygons of cubes in the cubical tessellation of $\E$ with vertex set $a\mathbb{Z}^3$, with the base face $F_2$ unchanged. Here each cube occurs and contributes all four Petrie polygons, forming again a finite subcomplex $\{6,3\}_4$. In particular, the edge graph of $\K_{6}(1,2)$ is the full edge graph of the cubical tessellation and every edge of $\K_{6}(1,2)$ lies in eight faces, two coming from each cube that contains it. The vertex-figure of $\K_{6}(1,2)$ at $o$ is the double-edge graph of the octahedron with vertices $(\pm a,0,0)$, $(0,\pm a,0)$, $(0,0,\pm a)$, and the vertex-figure group $\langle R_1,G_2\rangle$ is $[3,4]$.  

Note that the faces of $\mathcal K_6(1,2)$ are those of the Petrie duals of the facets $\{\infty, 3\}_6 \# \{\, \}$ of the $4$-apeirotope $\apeir\{3,4\}$, with $\{3,4\}$ properly chosen (see (\ref{4apeirotopes})).

\medskip
\noindent{\em Case IIb: $E_1'$ is perpendicular to a main diagonal of $C$}
\medskip

Recall that $E_2'$, $R_0'$ and $R_1$ are exactly as in the previous case. Since the rotation axis of $R_0'$ lies in $E_1'$, there is exactly one choice for $E_1'$ (as usual, up to congruence), namely the plane through $o$ perpendicular to the main diagonal connecting $\pm (1,-1,1)$. Hence $G_2$ is cyclic of order $3$ or dihedral of order $6$, and the twin vertex is of the form $v=(a,-a,a)$ with $a\neq 0$. 

We can immediately rule out the possibility that $G_2$ is cyclic. In fact, by Lemma \ref{casesnum1}, the presence of a $3$-fold rotation, $S$, in $G_2$ already forces $G_2$ to be dihedral. Observe here that the vertex-figure group $\langle R_1,G_2\rangle$ of $\K$ at $o$ is the subgroup $[3,3]$ of $[3,4]$ and hence contains the full dihedral subgroup $D_3$ of $[3,3]$ that contains $S$; moreover, this subgroup $D_3$ fixes $o$ and $v$, so that $G_2=D_3$. Thus $G_2$ cannot be cyclic.

\begin{figure}
\begin{center}
\includegraphics[width=9cm, height=7cm]{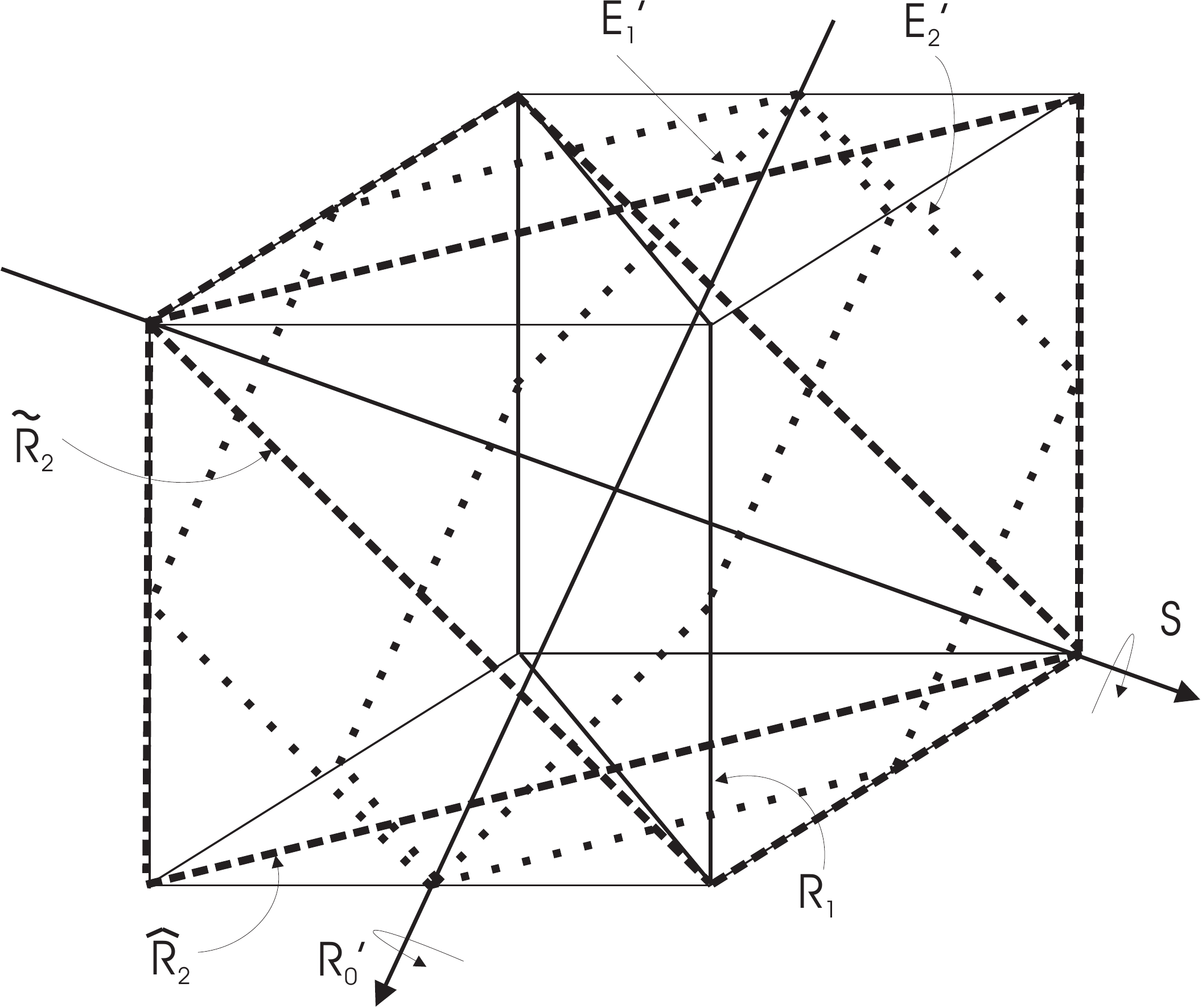}
\caption{The special group of the complex $\K_{7}(1, 2)$}\label{k712}
\end{center}
\end{figure}

Now suppose $G_2$ is dihedral of order $6$. Then we are necessarily in case (B). As generators for $G_2$ we may take the reflections $\hat{R}_{2}$ and ${\widetilde R}_2$ in the planes $y=-x$ and $z=-y$, respectively. This leads to the generators $R_0,R_1,{\hat R}_{2},{\widetilde R}_2$ of $G$ given by
\begin{equation}
\label{7gen12}
\begin{array}{rccl}
R_0\colon    & (x,y,z)  &\mapsto & (-z,-y,-x) + (a,-a,a),\\
R_{1}\colon & (x,y,z)  &\mapsto  & (y,x,z),\\
\hat{R}_{2}\colon  & (x,y,z) &\mapsto & (-y,-x,z),\\
\widetilde{R}_2\colon & (x,y,z)  &\mapsto & (x,-z,-y),
\end{array}
\end{equation}
with $a\neq 0$ (see Figure \ref{k712}). The resulting regular complex, denoted $\K_{7}(1,2)$, has skew hexagonal faces such that six surround each edge.

The geometry of $\K_{7}(1,2)$ is more complicated than in the previous cases. The vertex-set is 
$2a\Lambda_{(1, 1, 0)} \cup ((a, -a, a) + 2a\Lambda_{(1, 1, 0)})$. The edges of $\K_{7}(1,2)$ run along main diagonals of the cubes in the cubical tessellation with vertex-set $a\mathbb{Z}^3$, so in particular the faces are not Petrie polygons of cubes. The base face $F_2$ has vertices 
\[ \begin{array}{c}
v_0:=o,\;\, v_{1}:=v=(a,-a,a),\;\, v_{2}:=(0,-2a,2a),\; \\[.03in]
v_{3}:=(-a,-a,3a),\;\, v_{4}:=(-2a,0,2a),\;\, v_{5}:=(-a,a,a), 
\end{array} \]
in this order. The vertex-figure group $\langle R_1,G_2\rangle$ of $\K_{7}(1,2)$ at $o$ is again $[3,3]$. The vertex-figure at $o$ itself is the double-edge graph of the tetrahedron with vertices $(a,-a,a)$, $(-a,a,a)$, $(a,a,-a)$, $(-a,-a,-a)$. 

The twelve faces of $\K_{7}(1,2)$ containing the vertex $o$ can be visualized as follows. Consider the three (nested) cubes $aC$, $2aC$ and $3aC$, referred to as the {\em inner\/}, {\em middle\/} or {\em outer\/} cube, respectively; here $C$ is our reference cube with vertices $(\pm 1,\pm 1,\pm 1)$. Recall that there are two ways of inscribing a regular tetrahedron at alternating vertices of the inner cube, one given by the tetrahedron $P_a$ (say) with vertices $(a,-a,a)$, $(-a,a,a)$, $(a,a,-a)$, $(-a,-a,-a)$ that occurs as the vertex-figure at $o$. Clearly, $P_a$ is the convex hull of the orbit of $v_1$ under the vertex-figure group $[3,3]$. Let $P_{2a}$ and $P_{3a}$, respectively, denote the convex hulls of the orbits of $v_2$ and $v_3$ under $[3,3]$. Then $P_{2a}$ is the cuboctahedron whose vertices are the midpoints of the edges of the middle cube, and $P_{3a}$ is a truncated tetrahedron with its vertices on the boundary of the outer cube (and with four triangular and four hexagonal faces). 

With these reference figures in place, a typical hexagonal face of $\K_{7}(1,2)$ with vertex $o$ then takes its other vertices from $P_a$, $P_{2a}$ and $P_{3a}$, such that the vertices adjacent to $o$ are vertices of $P_a$, those two steps apart from $o$ are vertices of $P_{2a}$, and the vertex opposite to $o$ is a vertex of $P_{3a}$. In particular, for $F_2$, the vertices $v_1,v_5$ are from $P_a$, the vertices $v_2,v_4$ are from $P_{2a}$, and vertex $v_3$ is from $P_{3a}$. More precisely, since the vertex-figure is the double-edge graph of the tetrahedron, every edge $e$ of $P_a$ determines two faces of $\K_{7}(1,2)$ as follows. First note that the midpoint of $e$ lies on a coordinate axis and determines a positive or negative coordinate direction. Moving along this direction to the boundary of $P_{2a}$, we encounter a square face of $P_{2a}$ with one pair of opposite edges given by translates of $e$. Then each edge in this pair determines a unique face of $\K_{7}(1,2)$ that contains $o$. Thus there are two faces associated with $e$. To find their final vertex we move even further to the boundary of $P_{3a}$ until we hit an edge of $P_{3a}$. The vertices of this edge then are the opposites of $o$ in the two faces determined by $e$. Note here that the orbit of an edge of $P_{2a}$ under $[3,3]$ consists only of twelve edges, namely those of four mutually non-intersecting triangular faces of $P_{2a}$. On the other hand, $P_{3a}$ has exactly twelve vertices. 

Moreover, notice that the vertex-figure at a vertex adjacent to $o$ (such as $v$) similarly uses the other regular tetrahedron inscribed at alternating vertices of the inner cube. Thus both tetrahedra inscribed in the inner cube occur (in fact, already at every pair of adjacent vertices of $\K_{7}(1,2)$). Rephrased in a different way, every vertex $w$ of $\K_{7}(1,2)$ lies in four edges of $\K_{7}(1,2)$, such that their (outer) direction vectors point from $w$ to the vertices of a tetrahedron that is the translate by $w$ of a regular tetrahedron inscribed in the inner cube $aC$, and such that adjacent vertices of $\K_{7}(1,2)$ always use different inscribed tetrahedra. 

Finally, observe that the faces of $\mathcal K_7(1,2)$ are those of the Petrie duals of the facets $\{\infty, 3\}_6 \# \{\, \}$ of the $4$-apeirotope $\apeir\{3,3\}$, with $\{3,3\}$ properly chosen. The common edge graph of $\apeir\{3,3\}$ and $\mathcal K_7(1,2)$ is the famous {\em diamond net\/} (see \cite[p. 241]{rap}). The latter models the diamond crystal, with the carbon atoms sitting at the vertices and with the bonds between adjacent atoms represented by the edges (see also 
\cite[pp. 117,118]{wells}).

\medskip
\noindent{\em Case IIc: $E_1'$ is perpendicular to the line through the midpoints of a pair of antipodal edges of $C$}
\medskip

Since the rotation axis of $R_0'$ must lie in $E_1'$, there is exactly one choice for $E_1'$, namely the plane through $o$ perpendicular to the line through $\pm (1,0,1)$. Hence $G_2$ is cyclic of order $2$ or dihedral of order $4$, and the twin vertex is of the form $v=(a,0,a)$ with $a\neq 0$. 

If $G_2$ is cyclic of order $2$, then $\K$ is a regular polyhedron with complete dimension vector $(1,2,1)$. Its faces are skew hexagons, but now the edges are parallel to face diagonals of the cube; in fact, the faces of $\K$ are congruent to Petrie polygons of regular octahedra (see also $\K_{8}(1,2)$ described below). Moreover, $\K$ has planar hexagonal vertex-figures, since $R_1S$ is a rotatory reflection of period $6$ (with invariant plane perpendicular to the line through $\pm (1,1,-1)$). It follows that $\K$ must be the polyhedron $\{6,6\}_4$ (see \cite[p.225]{rap}).

\begin{figure}
\begin{center}
\includegraphics[width=8cm, height=7cm]{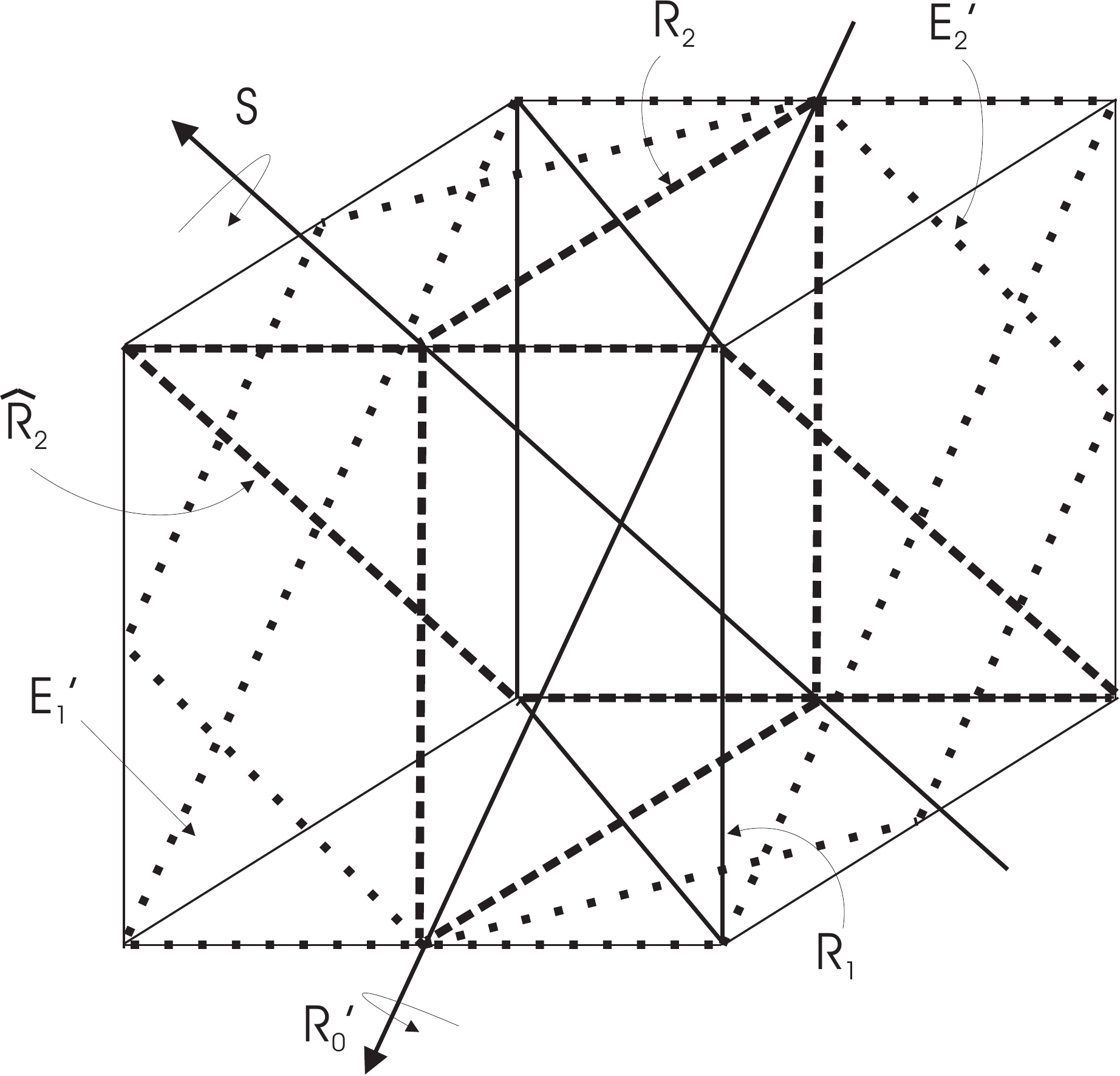}
\caption{The special group of the complex $\K_{8}(1, 2)$}\label{k812}
\end{center}
\end{figure}

Now suppose $G_2$ is dihedral of order $4$. Then $G_2$ is generated by the reflections $R_2$ and $\hat{R}_{2}$ in the planes $y=0$ and $z=x$, respectively, and we are in case (A). The generators $R_0,R_1,R_2,\hat{R}_{2}$ of $G$ are given by
\begin{equation}
\label{8gen12}
\begin{array}{rccl}
R_0\colon    & (x,y,z)  &\mapsto & (-z,-y,-x) + (a,0,a),\\
R_{1}\colon & (x,y,z)  &\mapsto  & (y,x,z),\\
R_{2}\colon  & (x,y,z) &\mapsto & (x,-y,z),\\
\hat{R}_2\colon & (x,y,z)  &\mapsto & (z,y,x),
\end{array}
\end{equation}
where again $a\neq 0$ (see Figure \ref{k812}). Now we obtain a regular complex, denoted $\K_{8}(1,2)$, with skew hexagonal faces, four surrounding each edge. This contains the polyhedron $\{6,6\}_4$ as a subcomplex.

The vertex-set of $\K_{8}(1,2)$ is $\Lambda_{(a,a,0)}$. Relative to the cubical tessellation with vertex-set $a\mathbb{Z}^3$, the edges of $\K_{8}(1,2)$ are face diagonals and the faces of $\K_{8}(1,2)$ are Petrie polygons of the octahedral vertex-figures of the tessellation at the vertices that are not in $\Lambda_{(a,a,0)}$.  At these octahedral vertex-figures, all Petrie polygons occur as faces of $\K_{8}(1,2)$, forming a finite subcomplex $\{6,4\}_3$ (the Petrie dual of $\{3,4\}$). For example, the base face $F_2$ has vertices $o$, $v=(a,0,a)$, $(0,-a,a)$, $(0,0,2a)$, $(-a,0,a)$, $(0,a,a)$, in this order, and is a Petrie polygon of the octahedral vertex-figure at $(0,0,a)$, its center. Thus every edge of $\K_{8}(1,2)$ is surrounded by four faces, arising in pairs from the octahedral vertex-figures at the two remaining vertices of the square face of the tessellation that has the edge as a face diagonal. Moreover, the vertex-figure group $\langle R_1,G_2\rangle$ of $\K_{8}(1,2)$ at $o$ is the octahedral group $[3,4]$, occurring here with (standard) generators $R_2,R_1,\hat{R}_2$. The vertex-figure of $\K_{8}(1,2)$ at $o$ is the (simple) edge graph of the cuboctahedron with vertices $(\pm a,\pm a, 0)$, $(0,\pm a,\pm a)$, $(\pm a,0,\pm a)$. 
\bigskip

Finally, then, this completes our investigation of all possible relative positions of mirrors of basic generators for the special group $G_{*}=[3,4]$ under the assumption that the mirror vector is $(1,2)$. In particular, we can summarize our results as follows.

\begin{theorem}\label{thm12}
Apart from polyhedra, the complexes $\K_{1}(1,2),\ldots,K_{8}(1,2)$ described in this section are the only simply flag-transitive regular polygonal complexes with finite faces and mirror vector $(1,2)$, up to similarity.
\end{theorem}

\section{Complexes with infinite faces and mirror vector~$(1, 2)$}
\label{infinmir12}

In this section we complete the enumeration of the simply flag-transitive regular complexes with  mirror vector $(1, 2)$. We saw in the previous section that, apart from polyhedra, there are only eight such complexes with finite faces. Now we prove that, again apart from polyhedra, there are no examples with infinite faces.

Let $\K$ be a simply flag-transitive regular complex with an affinely irreducible group $G = \langle R_0, R_1, G_2 \rangle$, where again $R_0$, $R_1$, $G_2$ are as in Lemma~\ref{r0r1}, the mirror vector is $(1,2)$, and the base vertex of $\K$ is $o$. Then the vertex-figure group $\langle R_1,G_2\rangle$ of $\K$ at $o$ is a subgroup of $G_*$, since both generators fix the base vertex, $o$. Recall that $v:=oR_0$ is called the twin vertex of $\K$.

Now suppose $\K$ has infinite faces. Then the faces must be planar zigzags. In fact, since the base face $F_2$ is an infinite polygon, the mirrors of the half-turn $R_0$ and plane reflection $R_1$ cannot meet and must be perpendicular to the plane $E$ spanned by $o$, $v$ and $vR_1$. Hence $F_2$ is a zigzag contained in $E$, and the mirrors of $R_0$ and $R_1$ are parallel. 

Now consider the elements $R_0'$ and $R_{1}$ of the special group $G_*$. Since the axis of $R_0'$ lies in the mirror of $R_1$, their product $R_{0}'R_{1}$ is the reflection in the plane through $o$ that is perpendicular to the mirror of $R_1$ and contains the axis of $R_0'$. Note here that $E$ is the plane through $o$ that is perpendicular to both mirrors, of $R_1$ and $R_{0}'R_{1}$. The twin vertex $v$ lies in $E$ but not on either mirror, since $v$ is not fixed by $R_1$ or $R_{0}'R_{1}$. Thus $G_*$ contains a pair of reflections with perpendicular mirrors and must be a group $[3,3]^*$, $[3,3]$ or $[3,4]$ (see (\ref{finthree})). Moreover, $G_2$ is cyclic or dihedral of order $r$. As for complexes with finite faces (and for exactly the same reasons), $G_2$ must contain a non-trivial rotation, $S$ (say), generating its rotation subgroup.  

Again we use the cube $C$ with vertices $(\pm 1,\pm 1,\pm 1)$ as the reference figure for the action of $G^*$. As in Case~I of the previous section, there are two possible choices for the half-turn $R_0'$:\ either $R_0'$ rotates about the center of a face of $C$ or about the midpoint of an edge of $C$. In the following we treat these as Case~I and Case~II. In any case, the rotation axis of $R_0'$ must lie in the mirrors of $R_1$ and $R_{0}'R_{1}$.

\medskip
\noindent{\em Case I:  $R_0'$ rotates about the center of a face of $C$}
\medskip

Suppose $R_0'$ rotates about the $y$-axis (say). Then there are two possible configurations of mirrors of $R_1$ and $R_{0}'R_{1}$:\  either they are the $xy$-plane and $yz$-plane, or the planes $z=x$ and $z=-x$. In either case, $E$ is the $xz$-plane plane and contains the rotation axis of $S$ (passing through $o$ and $v$). 

The first configuration can be ruled out immediately; in fact, then $S$ would rotate about the midpoint of an edge of $C$ (recall that $v$ cannot lie on either mirror) and $G_*$ would be a reducible group with invariant plane $E$, regardless of whether $G_2$ is cyclic or dihedral. This only leaves the second configuration, for which $S$ must also rotate about the center of a face of $C$. 

For the second configuration we can immediately exclude the possibility that $G_2$ is cyclic of order $2$ or dihedral of order $4$ with plane mirrors given by coordinate planes, once again appealing to irreducibility.  Furthermore, $G_2$ can also not be cyclic of order $4$ or dihedral of order $8$, as can be seen as follows.  First observe that, if $S$ has period $4$, then $G_2$ must actually be dihedral of order $8$. In fact, then the vertex-figure group $\langle R_1,S\rangle$ of $\K$ at $o$ is the full group $[3,4]$ and hence contains the full dihedral subgroup $D_4$ that contains $S$; but this subgroup $D_4$ also fixes $v$, so that it must coincide with $G_2$. On the other hand, if $G_2$ is dihedral of order $8$, then $G_2$ contains the reflection with mirror $E$ and hence $E$ is a face mirror of $\K$; in other words, the reflection in $E$ fixes the base flag of $\K$, in contradiction to our assumption that $\K$ is simply flag-transitive. Recall that the existence of face mirrors for $G$ immediately forces $G$ to be the full symmetry group of a regular $4$-apeirotope in $\E$ (see Lemma~\ref{reflectionface} and Theorem~\ref{nontriv}). Thus we have eliminated all but one possibility for $G_2$.

It remains to also reject the final possibility that $G_2$ is dihedral of order $4$ with plane mirrors determined by face diagonals of $C$. In this case the configuration of mirrors and rotation axes for the generators suggests that $\K$ would have to coincide with the $2$-skeleton of the regular $4$-apeirotope 
\[ {\cal P}:= \{\{\infty,3\}_6\# \{\, \},\{3,4\}\} = \apeir\{3,4\} \]
(see \cite[Ch.7F]{rap}). However, the $2$-skeleton is not a simply flag-transitive complex, so that the final case can also be excluded once this conjecture about $\K$ has been verified. We can accomplish the latter by employing Wythoff's construction as follows. Here we assume that the twin vertex is given by $v=(0,0,a)$ with $a\neq 0$.

Let $T_0,T_1,T_2,T_3$ denote distinguished generators for the symmetry group $G({\cal P})$ of $\cal P$, taken such that $T_1,T_2,T_3$ are distinguished generators of the octahedron $\{3,4\}$ with vertices at the face centers of the scaled reference cube $aC$, and such that $T_0$ is the point reflection in the point $\frac{1}{2}v$ halfway between the base vertex $o$ and twin vertex $v$ shared by $\K$ and $\cal P$. Here we may further assume that the generators $T_1,T_2,T_3$ are chosen such that $T_1=R_1$, $S=T_{2}T_{3}$, and $T_3$ is the plane reflection with mirror $E$ (the $xz$-plane). Then it is immediately clear that $T_3$ keeps the base face $F_2$ of $\K$ pointwise fixed, since the latter lies in $E$. Moreover, the vertex-figure group $\langle R_1,S\rangle$ of $\K$ at $o$ is $[3,3]$ and has index $2$ in the vertex-figure group $\langle T_1,T_2,T_3\rangle = [3,4]$ of $\cal P$ at $o$ (with $T_3$ representing the non-trivial coset), and $T_0=T_{3}R_{0}$. In fact, $G$ itself is a subgroup of $G({\cal P})$ of index at most $2$ (with $T_3$ representing the non-trivial coset if the index is $2$), since conjugation by $T_3$ leaves $G$ invariant. Now recall that the apeirotope $\cal P$ can be constructed by Wythoff's construction applied to $G({\cal P})$ with initial vertex $o$ (see \cite[Chs.5A,7F]{rap}); in particular, its base vertex is $o$, its base edge is $\{o,v\}$ (the orbit of $o$ under $\langle T_0\rangle$), its base $2$-face is the orbit of the base edge under $\langle T_0,T_1\rangle$, and its base $3$-face is the orbit of the base $2$-face under $\langle T_0,T_1,T_2\rangle$. The given complex $\K$ can similarly be derived by Wythoff's construction applied to $G=G(\K)$ with the same initial vertex, $o$.

We claim that $\K$ coincides with the $2$-skeleton of $\cal P$. First we verify that $\K$ and the $2$-skeleton of $\cal P$ have the same base flag. Since they obviously share their base vertices and base edges, we only need to examine their base $2$-faces; but clearly these must also agree, since $T_3$ fixes $E$ pointwise and $R_{0}=T_{3}T_{0}$. Now recall that the flags of a regular complex are just the images of the base flag under the respective group. If follows that $\K$ and the $2$-skeleton of $\cal P$ must indeed be identical, since $G$ is a subgroup of $G({\cal P})$ of index at most $2$ and the only nontrivial coset of $G$ (if any) is represented by an element, $T_3$, that fixes the base flag. (It also follows at this point that $G$ must indeed have index $2$ in $G({\cal P})$, since otherwise their vertex-figure subgroups would agree.)

\medskip
\noindent{\em Case II:  $R_0'$ rotates about the midpoint of an edge of $C$}
\medskip

Suppose $R_0'$ is the half-turn about the line through $o$ and $(1,1,0)$ (say). Then the mirrors of $R_1$ and $R_{0}'R_{1}$ are the $xy$-plane and the plane $y=x$, respectively, or vice versa. In either case, $E$ is the plane $y=-x$ and $S$ is a $3$-fold rotation about a main diagonal of $C$ contained in $E$. 

First we deal with the case that $R_1$ is the reflection in the plane $y=x$. Since $S$ is a $3$-fold rotation, the vertex-figure group $\langle R_1,S\rangle$ of $\K$ at $o$ is the full tetrahedral group $[3,3]$ and hence contains the full dihedral subgroup $D_3$ containing $S$. In particular, the reflection in $E$ belongs to $D_3$, so that $E$ is a face mirror for $\K$; however, this is impossible, since $\K$ is simply flag-transitive (see Lemma~\ref{reflectionface} and Theorem~\ref{nontriv}).

Next suppose that $R_1$ is the reflection in the $xy$-plane. We can immediately exclude the possibility that $G_2$ is dihedral of order $6$ (although, strictly speaking, this is subsumed under the next case). In fact, if $G_2$ was dihedral, the reflection in $E$ would already lie in $G_2$ and hence $E$ would again be a face mirror of $\K$.

This only leaves the possibility that $G_2$ is cyclic of order $3$. Here the configuration of mirrors and rotation axes for the generators suggests that $\K$ would have to be the $2$-skeleton of the regular $4$-apeirotope 
\[ {\cal P}:= \{\{\infty,4\}_4\# \{\, \},\{4,3\}\} = \apeir\{4,3\} \]
(see \cite[Ch.7F]{rap}). Hence, if confirmed, this possibility can be ruled out as well, as then $K$ could not be a simply flag-transitive complex. The conjecture that $\K$ coincides with the $2$-skeleton of ${\cal P}$ can be verified as for the special configuration described under Case~I.  Let $T_0,T_1,T_2,T_3$ denote the distinguished generators for $G({\cal P})$, where $T_0$ is the point reflection in the point $\frac{1}{2}v$ (with $v=(a,-a,a)$, say), and the distinguished generators $T_1,T_2,T_3$ of the cube $\{4,3\}$ (taken as $aC$, say) are such that $T_1=R_1$, $S=T_{2}T_{3}$, and $T_3$ is the plane reflection with mirror $E$ (the plane $y=-x$). Now $\langle R_1,S\rangle$ is the subgroup $[3,3]^*$ of index $2$ in $\langle T_1,T_2,T_3\rangle = [4,3]$ (and does not contain $T_3$), and again $T_0=T_{3}R_{0}$. Then we can argue as before. The group $G$ is of index at most $2$ in $G({\cal P})$, and the reflection $T_3$ represents the non-trivial coset of $G$ (if any) and fixes the common base flag of $\K$ and the $2$-skeleton of $\cal P$. Thus $\K$ coincides with this $2$-skeleton. (It also follows that $G$ must have index $2$ in $G({\cal P})$.)
\medskip

On a final note, the alert reader may be wondering why only two of the three possible $2$-skeletons of regular $4$-apeirotopes with infinite faces have occurred in our discussion (recall here Lemma~\ref{petrieop}). In fact, the $2$-skeleton of the third apeirotope $\{\{\infty,3\}_6\# \{\, \},\{3,3\}\}=\apeir\{3,3\}$ has occurred as well, at least implicitly, when we rejected the first possibility under Case~II that $R_1$ is the reflection in the plane $y=x$ (based on our observation that then $G$ would have to be the full group of this apeirotope since $E$ is a face mirror). 

In conclusion, we have established the following theorem.

\begin{theorem}\label{thm12infinite}
Apart from polyhedra, there are no simply flag-transitive regular polygonal complexes with infinite faces and mirror vector $(1,2)$.
\end{theorem}

The only regular polyhedra with infinite (zigzag) faces and mirror vector $(1,2)$ are the (blended) polyhedra $\{\infty,4\}_4 \# \{\,\}$, $\{\infty,6\}_3 \# \{\,\}$ and $\{\infty,3\}_6 \# \{\,\}$, obtained by blending the Petrie duals $\{\infty,4\}_4$, $\{\infty,6\}_3$ and $\{\infty,3\}_6$ of the regular plane tessellations $\{4,4\}$, $\{3,6\}$ and $\{6,3\}$, respectively, with the line segment $\{\,\}$ (see \cite[Ch.7E]{rap}). Each polyhedron is isomorphic to the corresponding Petrie dual (its plane component) but has its two sets of alternating vertices lying in two distinct parallel planes.
\vskip1in
\noindent
{\bf \Large Acknowledgment}\\

\noindent
We are very grateful to Peter McMullen for his comments on an earlier draft of this manuscript and for a number of very helpful suggestions. We would also like to thank an anonymous referee for a very thoughtful review with valuable suggestions for improvement.

\end{document}